\documentclass[11pt]{article}   % For Latex2e
\usepackage{amssymb,amscd,latexsym}   % For Latex2e
\usepackage{amsmath}
\usepackage{amsthm}
\usepackage{mathdots}
\usepackage[pagebackref]{hyperref}
\usepackage{color}
\usepackage{xcolor}
\newcommand{\rar}{\rightarrow}
\newcommand{\lar}{\longrightarrow}
\newcommand{\llar}{-\kern-5pt-\kern-5pt\longrightarrow}
\newcommand{\surjects}{\twoheadrightarrow}

\newtheorem{Theorem}{Theorem}[section]
\newtheorem{Lemma}[Theorem]{Lemma}
\newtheorem{Corollary}[Theorem]{Corollary}
\newtheorem{Proposition}[Theorem]{Proposition}
\newtheorem{Remark}[Theorem]{Remark}

\newtheorem{Conjecture}[Theorem]{Conjecture}

\newtheorem{Question}[Theorem]{Question}
\def\sqr#1#2{{\vcenter{\hrule height.#2pt
        \hbox{\vrule width.#2pt height#1pt \kern#1pt
            \vrule width.#2pt}
        \hrule height.#2pt}}}
\def\phi{\varphi}
\def\demo{\noindent{\bf Proof. }}
\def\square{\mathchoice\sqr64\sqr64\sqr{4}3\sqr{3}3}
\def\qed{\hspace*{\fill} $\square$}

\def\xx{{\bf x}}

\def\tt{{\bf t}}

\def\ZZ{{\bf Z}}
\def\ff{{\bf f}}

\def\ff{{\bf f}}

\def\aa{{\bf a}}
\def\qq{{\bf q}}
\def\vv{{\bf v}}

\def\codim{{\rm codim}\,}

\def\restr{{\kern-1pt\restriction\kern-1pt}}

\def\pp{{\mathbb P}}

\usepackage{tikz}
\usetikzlibrary{arrows,matrix,positioning}
\begin{document}
%\begin{titlepage}
\begin{center}
{\Large{\bf\sc Coordinate sections of generic Hankel matrices}}
\footnotetext{AMS 2010 Mathematics
Subject Classification (2010   Revision). Primary 13C40, 13D02, 13H10, 14E05, 14E07; Secondary 13C15, 14M10, 14M12, 14M15.}
\footnotetext{	{\em Key Words and Phrases}: Hankel matrix, Hessian, ideal of minors, polar map, gradient ideal, linear rank, special fiber, Cohen--Macaulay.}

\vspace{0.3in}

{\sc Rainelly Cunha}\footnote{Partially supported by a CAPES-PNPD post-doctoral fellowship (1723357/2017) from the Federal University of Sergipe, Brazil, during parts of this work. This author thanks the Department of Mathematics for providing the appropriate environment to carry the research.}
 \quad\quad
 {\sc Maral Mostafazadehfard}\footnote{Part of the work was done  while this author held a CAPES-PNPD Post-Doctoral fellowship ( 88882.317200/2019-01) from IMPA, Brazil. She thanks IMPA for providing the appropriate environment to carry the research.}
 \quad\quad
 {\sc Zaqueu Ramos}\footnote{Part of the work was done  while this author held  a grant from CAPES/FAPITEC (88887.157386/2017-00).}
 \quad\quad
 {\sc Aron  Simis}\footnote{Partially supported by a CNPq grant (302298/2014-2). Part of the work was done while this author held a PVNS Fellowship from CAPES (5742201241/2016), and subsequently  a Senior Visiting Research Grant at the ICMC-USP (S\~ao Carlos, SP). He thanks the Department of Mathematics of the Federal University of Paraiba and the ICMC-USP for providing an excellent environment for discussions on this work.}

\end{center}
%\end{titlepage}

%\tableofcontents
%\newpage

\begin{abstract}
	One deals with  degenerations by coordinate sections of the square generic Hankel  matrix over a field $k$ of characteristic zero, along with its main related structures, such as the determinant of the matrix, the ideal generated by its partial derivatives, the polar map defined by these derivatives, the Hessian matrix and the ideal of the submaximal minors of the matrix.  It is proved that the polar map is dominant for any such degenerations, and not homaloidal in the generic case. The problem of whether the determinant $f$ of the matrix is a factor of the Hessian with the (Segre) expected multiplicity is considered, for which the expected lower bound of the dual variety of $V(f)$ is established.
\end{abstract}

\section*{Introduction}

A good amount of algebraic work has addressed various aspects of Hankel matrices, if not to mention the non-enumerable list of papers dealing with its applications -- often under the disguise of  Toelitz matrices -- ranging from Functional Analysis to Orthogonal Polynomial Theory to the Moment Problem and Probability.
In commutative algebra and algebraic geometry they appear with functional entries, or actually polynomial entries.
An expressive case is that where the entries are actually linear forms, in fact even variables, giving room to some of the so-called determinantal ideals or varieties.
Typical geometric objects defined by Hankel matrices and their ideals of minors are secant varieties and normal scrolls.
Those in turn have many applications in discrete geometry, including graph theory.

In commutative algebra, the work of Watanabe (\cite{WatanabeJ}), Conca (\cite{Conca1998}), Eisenbud (\cite{Eisenbud1} and \cite{Eisenbud2}), among others, has a distinct role.
It is rarely the case that a paper on the subject does not quote one of these sources. 

The present work is a step into a theory of certain degenerations of the generic Hankel matrix. It is naturally inspired by previous work such as \cite{MAron} and \cite{Degen_generic}.
The first of these sources dealt with the ideal theoretic and homological nature of the generic square Hankel determinant and of its associated objects.
In addition there was a good deal of results on the so-called sub-Hankel matrix, a degeneration of the generic Hankel matrix, largely considered in \cite{CRS} for its significance in determinantal Cremona theory.

The second source above deals with certain degenerations of the square generic matrix. Except for the central idea of looking at degenerations of a generic object, the two have very little in common in regard to the expected results. 
Yet, a common thread in each was the goal to understand the properties of the  gradient (Jacobian) ideals, the corresponding Hessian determinants and the ideals of lower minors --  more particularly, the submaximal minors.
This thread will be equally pursued in this work.

Thus, the main theme of this work is the rough material coming from the first of the two sources, while keeping an eye on the fabrics of the second source.
That is, we focus on square Hankel matrices which are degenerations of the generic Hankel matrix, by a process of linear space sections to put in a more geometric way.

It may be more enlightening to start directly with a description of the contents of each section.

The first section deals with some preliminaries and explains the title of this work.
Roughly, given a generic member within some class of matrices, it consists in replacing each variable in a certain segment by zero. 
Lacking better names, we called the resulting matrix a degeneration by {\em coordinate sections} honoring \cite{GiMer}. The geometric idea behind this terminology is that one is looking at coordinate hyperplane sections of the associated determinantal varieties. 
This sort of degeneration has been thoroughly dealt with in \cite{Degen_generic} and \cite{Degen_sym}, and considered before by other authors (\cite{GiMer}, \cite{Eisenbud2}).
The main subject is thus a square Hankel matrix $\mathcal{H}_m[r]$ of order $m$, with $r$ zeros along the lower stretch of the its last column.
This section is about the basics of such matrices and the gradient ideals of the respective determinants.

Section~\ref{ideals_of_minors} deals with properties of the ideals of lower minors of the above degenerations.
We are particularly interested in the ideal of submaximal minors.
Several properties relating $I_{m-1}(\mathcal{H}_m)$ (here $\mathcal{H}_m=\mathcal{H}_m[0]$, the generic case) to the gradient ideal $J$ generated by partial derivatives of $\det \mathcal{H}_m$ are proved in \cite{MAron}. One of these is that $I_{m-1}(\mathcal{H}_m)$ is the minimal component of the primary decomposition of $J$. 
The main subject in the section is the codimension, primeness and Cohen--Macaulayness of ideals of minors of the degeneration $\mathcal{H}_m[r]$.
A central question addressed, but unfortunately still pending, is the precise structure of the ideal of submaximal minors and its associated algebras, such as the special fiber.

Having dealt for a while with ideals of minors of the degeneration matrix $\mathcal{H}_m[r]$, in Section~\ref{Hessian} we consider more closely its differential properties. Throughout and henceforth, we assume that the ground characteristic is zero. Thus, we look at the ideal $J=J(f)\in R$ generated by the partial derivatives of $f:=\det \mathcal{H}_m[r]$ -- usually called the {\em Jacobian ideal} of $f$, a terminology here deferred in favor of {\em gradient ideal} of $f$.
The Jacobian matrix of the partial derivatives is the so-called {\em Hessian matrix} of $f$.
By abuse, its determinant is called {\em the Hessian} of $f$.
First one considers the Hessian $h(f)$ of $f=\det \mathcal H{_m}[r]$.
Since we are assuming characteristic zero, the algebraic translation of its non-vanishing is that the analytic spread of the gradient ideal of $f$ is maximal, i.e., equal to $\dim R$ in this case.
The geometric impact of this result is that the polar map of $f$ is dominant for any value of $r$. It immediately raises the supplementary question as to whether the polar map is actually a Cremona map (i.e., whether $f$ is a homaloidal polynomial) -- this problem will be considered in Section~\ref{homaloidal}.
This state of affairs is quite surprising as the analogous degeneration in the case of the generic square matrix has a vanishing Hessian determinant as soon as $r\geq 1$ (\cite{Degen_generic}).

We observe that the extreme case $m-r=2$ (sub-Hankel) has been proved in \cite{CRS} by showing that the corresponding Hessian is a power of $x_{m+1}$ up to a nonzero coefficient from $k$.
At the other end of the spectrum, the fully generic case (i.e., $r=0$) has been given in \cite[Proposition 3.3.11]{Maral} by a method of degenerating the Hessian matrix by setting every variable to zero, except $x_m$.
Unfortunately, this simple degeneration does not work well all the way for arbitrary $r\geq 1$.
However, the idea of suitably degenerating the Hessian matrix can still be used as is shown in the Appendix.
The proof  draws on degenerating the Hessian matrix by as many  coordinate sections as possible. It is not apparent what is a best choice (if ever) since it may depend on both $m$ and $r$. Some choices may work in one case, but not elsewhere. 
For any promising degeneration the argument will be fatally long.
Luckily, in the case $r=0$ (i.e., the fully generic case), one can alternatively look at the more elegant argument in Corollary~\ref{Hessian_is_not_null}.

An interesting question in general is whether $f:=\det {\mathcal H}_m[r]$ is a factor of its Hessian determinant $h(f)$ with multiplicity $\geq 1$. 
If this is the case, then $f$ is said in addition to have the {\em expected multiplicity} (according to Segre) if its multiplicity as a factor of $h(f)$ is  $ {\rm codim}(V(f)^*)-1$, where $V(f)^*$ denotes the dual variety to the hypersurface $V(f)$ (see \cite{CRS}).
In a similar context this question has been tackled in \cite[Section 3.2]{MAron} (see also \cite[Remark 2.7 (2)]{Degen_generic}) for a couple of corrections). In the case of the $m\times m$ generic matrix $\mathcal G$ the Hessian of $f=\det \mathcal G$ equals $f^{m(m-2)}$. Now, it is well-known that the dual variety to $f$ in this case is the ideal of $2$-minors of the same matrix in the dual variables.
Since the latter has codimension $(m-1)^2=m(m-2)+1$ then $f$ is a divisor of its Hessian determinant with the expected multiplicity.
An extension of this result has been conjectured in \cite{MAron}, namely for the class of higher leap Hankel-wise $m\times m$ matrices, and confirmed for a few low values of $m$ and the leap. 

In a different north, the question comes up for degenerations of $\mathcal G$ that preserve the non-vanishing of the corresponding Hessian determinant.
A case has been dealt with in \cite[Section 2.3]{Degen_generic}, where the same phenomenon has been proved to hold.
On the other hand, for the analogue of the  degeneration studied here as applied to $\mathcal G$ the Hessian always vanishes for $r\geq 1$, so the question does not come up.
It is somewhat surprising that the degenerations $\mathcal H_m[r]$ by coordinate sections recover the non-vanishing of the corresponding Hessian determinant as mentioned above. We prove that $\dim V(f)^* $ is at least $m-1$ and state some conjectured assertions, followed by some justification.

The second topic of Section~\ref{Hessian} is  the gradient ideal $J$.
Here the main result is the structure of its minimal prime ideals in the case where $1\leq r\leq m-3$. This question for the cases where $r=0$ (generic) or $r=m-2$ has been settled in \cite{MAron}.
In contrast to the generic case, $J$ is never a reduction of the ideal of submaximal minors. This lies at the crux of the typical difficulty in handling the case of $1\leq r\leq m-3$.
The complete set of associated primes of $J$ remains a mystery even in the generic case. 

The last topic of the section is the linear behavior of the gradient ideal $J$.
The two notions of linear behavior for a homogeneous ideal -- its linear rank and its linear type property -- have no apparent relation in general. Yet, for a homogeneous ideal generated by its linear system on the initial degree, it has been noted that the two notions intertwine as regards the birationality of the map defined by this linear system  (see \cite[Theorem 3.2 and Proposition 3.4]{AHA} and earlier references mentioned there).

For the generic Hankel matrix $\mathcal{H}_m$, at least in null characteristic,  it was proved that the linear rank of $J$ is $3$ (\cite[Theorem 3.3.5]{Maral}).
At the other end of the spectrum, for the degeneration $\mathcal{H}_m[m-2]$ it was proved that the linear rank of $J$ is maximal possible ($=m$) (\cite{CRS}) and $J$ is in addition of linear type (\cite[Section 4.1, Theorem 4.8]{MAron}).
Here we conjecture that, for $1\leq r\leq m-3$, the linear rank of $J$ is $2$.
We give a proof of this assertion based on another conjecture regarding regular sequence modulo the gradient ideal in the generic case.
The linear rank conjecture itself fails in positive characteristic. 

Section~\ref{generic_case_revisited} is mainly about properties of the generic Hankel matrix, as a natural sequel to questions posed in \cite{MAron}.
Thus, we first make a thorough review of the background combinatorics, including Pl\"ucker relations, that affects this case.
The main results are Proposition~\ref{Hankel_is_Grassmann} and Theorem~\ref{non-homaloidal_main}. The first essentially states that the defining equations of the algebra of submaximal minors are the same as those of a Grassmannian (namely, the Pl\"ucker relations) -- a somewhat surprising result.
The second tells that the polar map in the generic case is not homaloidal and, in addition, that the gradient ideal is a reduction of the ideal of submaximal minors.
The first of these results has also been obtained by N. Medeiros by a geometric argument.

The last part of the section contains a discussion of the elements to prove a parallel result to the effect that also in the case $1\leq r\leq m-3$ the polar map is not homaloidal either -- thus leaving us with the sole case where the polar map is homaloidal, namely, $r=m-2$ (cf. \cite{CRS}).

\section{Preliminaries}

\subsection{Recap of basic notions}

Let $(R,\mathfrak{m})$ denote a Notherian local ring and its maximal ideal (respectively, a standard graded ring over a field and its irrelevant ideal).
For an ideal $I\subset \mathfrak{m}$ (respectively, a homogeneous ideal $I\subset \mathfrak{m}$), the \emph{special fiber} of $I$ is the ring $\mathcal{R}(I)/\mathfrak{m}\mathcal{R}(I)$.
Note that this is an algebra over the residue field of $R$.
The (Krull) dimension of this algebra is called the {\em analytic spread} of $I$ and is denoted $\ell (I)$. 

Quite generally, given ideals $J\subset I$ be ideals in a ring $R$,  $J$ is said to be a \emph{reduction} of $I$ if there exists an integer $n\geq 0$ such that $I^{n+1}=JI^n.$
A reduction $J$ of $I$ is called \emph{minimal} if no ideal strictly contained in $J$ is a reduction of $I$.
The \emph{reduction number} of $I$ with respect to a reduction $J$ is the minimum integer $n$ such that $JI^{n}=I^{n+1}$. It is denoted by $\mathrm{red}_{J}(I)$. The (absolute) \emph{reduction number} of $I$ is defined as $\mathrm{red}(I)=\mathrm{min}\{\mathrm{red}_{J}(I)~|~J\subset I~\mathrm{is}~\mathrm{a}~\mathrm{minimal}~\mathrm{reduction}~\mathrm{of}~I\}.$ 
If $R/\mathfrak{m}$ is infinite every minimal reduction of $I$ is minimally generated by exactly $\ell(I)$ elements. In particular, every reduction of $I$ contains a reduction generated by $\ell(I)$ elements.
The following invariants are related in the case of $(R,\mathfrak{m})$:
$$\mathrm{codim}(I)\leq \ell(I) \leq \min\{\mu(I), \mathrm{dim}(R)\},$$
where $\mu(I)$ stands for the minimal number of generators of $I$.
If the rightmost inequality turns out to be an equality, one says that $I$ has maximal analytic spread.
By and large, the ideals considered in this work will have $\dim R\leq \mu(I)$, hence being of maximal analytic spread means in this case that $\ell(I)=\dim R$.

Suppose now that  $R$ is standard graded over a field $k$ and $I$ is an ideal of grade $\geq 1$ generated by $n+1$ forms of a given degree $s$. One has a free graded presentation $$ R(-(s+1))^t \; \oplus\; \sum_{j\geq 2} R(-(s+j))\stackrel{\varphi}{\longrightarrow} R(-s)^{n+1}\longrightarrow I\longrightarrow 0  $$
for suitable shifts $-(s+j)$ and $t\geq 0$.   Of much interest is the image of $R(-(s+1))^t$ by $\varphi$, so-called {\it  linear part of $\varphi$} -- often denoted $\varphi_1$.   Since $\phi$ has a rank, so does $\phi_1$. One says that  the rank of $\varphi_1$ is the {\it linear rank of} $\varphi$ (or of $I$)  and that $\varphi$ has {\it maximal linear  rank} provided its linear rank is $n$ (=rank($\varphi$)). Clearly, the latter condition  is trivially satisfied if $\varphi=\varphi_1 $, in which case $I$ is said  to have {\it linear presentation} (or is {\it linearly presented}).

 Given a monomial order the polynomial ring $R$ over a field, if $f\in R$ we denote by ${\rm in}(f)$ the initial term of $f$  and by ${\rm in}(I)$  the ideal generated by the initial terms of the elements of $I$, called the initial ideal of $I$. 
For the general theory of monomial ideals and Gr\"obner bases we refer to \cite{HeHiBook}.

\subsection{Hankel matrix tools}

The generic  Hankel matrix of order $s\times (n-s+1)$ in $n$ variables is the catalecticant 
\begin{equation}\label{Hankel_nonsquare}
\mathcal{H}_{s, n-s+1}:=\left(
\begin{matrix}
x_1&x_2&x_3&\ldots &x_s&\ldots&x_{n-s}&x_{n-s+1}\\
x_2&x_3&x_4&\ldots &x_{s+1}&\ldots&x_{n-s+1}&x_{n-s+2}\\
x_3&x_4&x_5&\ldots &x_{s+2}&\ldots &x_{n-s+2}&x_{n-s+3}\\
\vdots &\vdots &\vdots&\ldots &\vdots &\vdots&\vdots \\
x_{s}&x_{s+1}&x_{s+2}&\ldots &x_{2s-1}&\ldots &x_{n-1}&x_{n}\\
\end{matrix}
\right),
\end{equation}
where, say,  $s\leq n-s+1$ (i.e., $n\geq 2s-1$). 

Specially notable is the case where $s=n-s+1:=m$, the square Hankel matrix of order $m$, henceforth denotes $\mathcal{H}_m$.
This work focuses on certain degenerations of the generic  Hankel matrix.

As usual, the ideal of $t$-minors of an arbitrary matrix $\mathcal{M}$ will be denoted $I_t(\mathcal{M})$.

The ideals of minors of a Hankel matrix have a notable behavior.
It has been proved in \cite[Proposition 4.3]{Eisenbud2}  that, for any $1\leq t\leq \min\{s,n-s+1\}$, the ideal of $t$-minors of  the generic Hankel matrix $\mathcal{H}_{s, n-s+1}$  is prime and of codimension $n-2t+2$. The proof of this result uses an important property of Hankel matrix  first made explicit in  the work of Gruson and Peskine:

\begin{Proposition}\label{GP} {\rm (\cite{GP})} For any $t\leq s$, one has
$I_t(\mathcal{H}_{s,n-s+1})=I_t(\mathcal{H}_{t,n-t+1})$.
\end{Proposition}

This property allows to reduce to the case of maximal minors.
In this case, to get the codimension one may use the fact that the Hankel matrix specializes to a well-known shape involving only $2m-2t+1$ variables -- but see Lemma~\ref{basic_minors} later on for a more precise statement. 
In particular, in the case of the square Hankel matrix $\mathcal{H}_m$, its ideal of submaximal minors coincides with the ideal of maximal minors of 
$\mathcal{H}_{m-1,,m+1}$. 

\begin{Remark}\label{same_vars}\rm It is perhaps worth observing that the notation of the Hankel matrix in (\ref{Hankel_nonsquare}) should include the set of variables used. Instead it is customary to fix the set of variables once for all, while administering various changes of sizes.
In this regard, given $\mathcal{H}_{s,n-s+1}$ and $t\leq s$ as in Proposition~\ref{GP}, a matrix such as $\mathcal{H}_{t,n-t+1}$ is uniquely defined.
\end{Remark}

The following result originally appeared in \cite{Golberg} in a different context. It has independently been obtained in \cite[Proposition 5.3.1]{Maral} in the presently stated form.

\begin{Proposition}\label{partials_vs_cofactors}
Let $\mathcal{M}$ denote a square matrix over $R=k[x_1,\ldots,x_n]$  such that  every entry is either 0 or $x_i$ for some $i=1,\ldots,n$. 
Then, for each $i=1,\ldots, n$, the partial derivative of  $f=\det(\mathcal{M})$ with respect to $x_i$ is the sum of the {\rm (}signed{\rm )} cofactors of the entry $x_i$ in all its slots as an entry of $\mathcal{M}$.
\end{Proposition}

\subsection{Degeneration by coordinate sections}\label{degens}

Let $R=k[\xx]=k[x_1,\ldots,x_n]$ denote a standard graded polynomial ring over a field $k$ of characteristic zero.
We will only consider degenerations induced by a $k$-algebra homomorphism of $R$ and, in fact, those induced by coordinate sections.
More particularly, our focus is on the following degenerations of the square version of (\ref{Hankel_nonsquare}):
$$
\left(
\begin{matrix}
x_1&x_2&\ldots &x_{m-1} & x_{m}\\
x_2&x_3&\ldots &x_{m}& x_{m+1}\\
\vdots &\vdots &\ldots &\vdots \\
x_{m-1}&x_{m}&\ldots &x_{2m-3} &x_{2m-2}\\
x_{m}&x_{m+1}&\ldots &x_{2m-2}& 0\\
\end{matrix}
\right),
\;
\left(
\begin{matrix}
x_1&x_2&\ldots &x_{m-2} &x_{m-1} & x_{m}\\
x_2&x_3&\ldots& x_{m-1} &x_{m}& x_{m+1}\\
\vdots &\vdots &\ldots &\ldots &\vdots &\vdots \\
x_{m-2}&x_{m-1}&\ldots &x_{2m-5} &x_{2m-4} &x_{2m-3}\\
x_{m-1}&x_{m}&\ldots &x_{2m-4} &x_{2m-3} &0\\
x_{m}&x_{m+1}&\ldots &x_{2m-3} &0& 0\\
\end{matrix}
\right),
$$
{\small
	$$
	\left(
	\begin{matrix}
	x_1&x_2&x_3&\ldots &x_{m-1} & x_{m}\\
	x_2&x_3& x_4& \ldots &x_{m}& x_{m+1}\\
	x_3&x_4& x_5 &\ldots  &x_{m+1}& 0\\
	\vdots &\vdots & \vdots  &\ldots &\vdots &\vdots \\
	x_{m-2}&x_{m-1}& x_m &\ldots  & 0 &0\\
	x_{m-1}&x_{m}& x_{m+1} & \ldots  & 0 &0\\
	x_{m}&x_{m+1}& 0 & \ldots  &0& 0\\
	\end{matrix}
	\right)
	$$}

We will denote by $\mathcal{H}_m[r]$ a Hankel degeneration as above, where $r$ denotes the number of zeros on the last column.
The last matrix in the above thread was dubbed {\em sub-Hankel} in \cite{CRS} (see also \cite{Maral}, \cite{MAron}).
This notation will also be used when the Hankel matrix is not necessarily square, namely,  $\mathcal{H}_{s,n-s+1}[r]$.
%However, in the non-square case we may need to map to zero a longer segment.

For a given $r$, the ground ring for $\mathcal{H}_m[r]$ is the polynomial ring $k[x_1,\ldots, x_{2m-r-1}]$.
If no confusion arises, when $r$ is fixed in the discussion, we will denote this ring simply by $R$.

It is quite clear that taking minors commute with ring homomorphisms.
In the present situation it is expressed in the following convenient way:

\begin{Lemma}\label{homo1}
Given an integer $r\geq 1$, let  $\Phi$ be the endomorphism of $R=k[x_1,\ldots,x_n]$ such that
$$\Phi(x_i)=\left\{
\begin{array}{ll}
0 & \mbox{if $i\geq 2m-2r$}\\
x_i & \mbox{otherwise}
\end{array}
\right.
$$	
Then
\begin{enumerate}
	\item[{\rm (a)}] $I_t(\mathcal{H}_{s,n-s+1}[r])=\Phi(I_t(\mathcal{H}_{s,n-s+1}))$, for all $1\leq t \leq s$.
	\item[{\rm (b)}] $I_t(\mathcal{H}_{s,n-s+1}[r])=I_t(\mathcal{H}_{t,n-t+1}[r])$, for all $1\leq t \leq s$.
\end{enumerate}
\end{Lemma}
\demo (a) This is clear.

(b) It follows from item (a) and  Proposition~\ref{GP}.
\qed

\medskip

\section{Ideals of minors}\label{ideals_of_minors}

\subsection{The determinant}

The proof  of the following proposition is inspired by an elementary fact observed in the case of the sub-Hankel in \cite[Remark 4.6 (c)]{CRS}, sufficiently generalized to the general case of a Hankel degeneration.

Actually, the observation will work for the generic Hankel matrix itself, thus avoiding drawing upon the general result about this matrix being $1$-generic \cite{Eisenbud2}.

\begin{Proposition}\label{the_determinant}
Let $\mathcal{H}_m[r]$ denote the  degeneration of the  $m\times m$ generic Hankel matrix considered in the previous section.
Let $R$ denote  the polynomial ring on the distinct non zero entries of the matrix. Then 
\begin{enumerate}
	\item[{\rm (i)}] $\det \mathcal{H}_m[r]\neq 0$.
	\item[{\rm (ii)}] 	$\det \mathcal{H}_m[r]\in R$ is irreducible if and only if $ r\leq m-2$.
\end{enumerate}
\end{Proposition}
\demo
Set $f:=\det \mathcal{H}_m[r]$.

(i) There are many elementary ways of verifying the non-vanishing of $f$.
Perhaps an easy one is to see that $f$ has a unique nonzero pure term in $x_m$, namely, the product of the entries along the main anti-diagonal.

(ii) The ``only if'' part is obvious since the determinant would then be a power of $x_m$ or zero.

For the reverse implication we will induct on $m$.
The initial step of the induction will be subsumed in the general step.
We may assume that $r\leq m-3$ since the case where $r=m-2$ has been established in \cite[Remark 4.6 (c)]{CRS}.

Since $x_1$ only appears once and on the first row, one easily sees that $f=x_1f_1+g$, for some $g\in k[x_2,\ldots,x_{2m-1-r}]$, where $f_1\in k[x_2,\ldots,x_{2m-1-r}]$ is the determinant of the Hankel degeneration of type $\mathcal{H}_{m-1}[r]$ obtained by omitting the first row and the first column of the original Hankel degeneration.

To show that $f$ is irreducible it suffices to prove that it is a primitive polynomial (of degree $1$) in $k[x_2,\ldots,x_{2m-1-r}][x_1]$.
Now, on one hand, $f_1$ is irreducible by the inductive hypothesis since one is assuming that $r\leq m-3=(m-1)-2$.
Therefore, it is enough to see that $f_1$ is not a factor of $g$.
For this, one verifies their initial terms in the revlex monomial order: ${\rm in}(f_1)=x_{m+1}^{m-1}$ and ${\rm in}(g)={\rm in}(f)=x_{m}^{m}$.
\qed

\begin{Remark}\rm
Since $f$ is homogeneous, an alternative argument for the case $r\leq m-3$ consists in showing that $R/(f)$ is normal. Since $R/(f)$ is a hypersurface ring, it suffices to prove that it is locally regular in codimension one. By Proposition~\ref{cod_of_gradient} below, proved independently, the gradient ideal $J$ has codimension $3=1+2$ provided $r\leq m-3$.
This proves that $f$ is irreducible when $r\leq m-3$.	
\end{Remark}

\subsection{Ideals of lower minors}

As a preliminary we need a result about the lower minors of a not necessarily square generic Hankel matrix.

Let $\mathcal H$ denote the generic Hankel matrix of order $a\times b$, with $a\leq b$ and entries $x_1,\ldots, x_{a+b-1}$. Set $R=k[x_1,\ldots,x_{a+b-1}]$ for the ground polynomial ring on these entries. Then the ideal $I_a(\mathcal H)$ of maximal minors has the expected codimension $b-a+1;$ in particular, $R/I_{a}(H)$ is a Cohen-Macaulay ring of dimension $2(a-1)$.
The expected codimension is a well-known result (see, e.g., \cite{Eisenbud1}).
Therefore, the Eagon--Northcoth principle is applicable to deduce the Cohen--Macaulayness and the dimension.

\begin{Lemma}\label{basic_minors} With the above notation, an explicit system of parameters of $R/I_{a}(H)$ is given by the residues of the elements $x_1,\ldots, x_{a-1}, x_{b+1},\ldots, x_{a+b-1}$.
\end{Lemma}
\demo 
Consider the corresponding fully generic $a\times b$ matrix $\mathcal{G}=(y_{i,j})$, for which it is well-known that $I_a(\mathcal{G})\subset S:=k[y_{i,j}]$ has codimension $b-a+1$. 
Let $D$ denote the ideal of $S$ generated by the $(a-1)(b-1)$ independent linear forms $y_{i,j}-y_{i+1,j-1}$, for $1\leq i< a, 2\leq j\leq b$.
Clearly, there is an isomorphism $S/D\simeq R$ specializing $\mathcal{G}$ to $\mathcal{H}$ and inducing an isomorphism
$$S/I_a(\mathcal{G})\, (\bmod D) \simeq R/I_a(\mathcal{H}).$$
Since the dimensions of $S/I_a(\mathcal{G})$ and $R/I_a(\mathcal{H})$ are $ab-(b-a+1)=(b+1)(a-1)$ and $a+b-1-(b-a+1)=2(a-1)$, respectively -- hence, their difference is precisely $(a-1)(b-1)$) -- and both rings are Cohen--Macaulay, it follows that $D$ is also a regular sequence on $S/I_a(\mathcal{G})$.
On the other hand, the initial ideal of $I_a(\mathcal{G})$ in the reverse lex order is generated by the products of the entries along the anti-diagonals of the $a$-minors. Therefore, it follows that the stated set $\{x_1,\ldots, x_{a-1}, x_{b+1},\ldots, x_{a+b-1}\}$ is a system of parameters (regular sequence) on $R/I_a(\mathcal{H})$.
\qed
\smallskip

For the next result we set $R:=k[x_1,\ldots, x_{2m-1}]$ and  $\bar{R}:=k[x_1,\ldots, x_{2m-r-1}]$, the latter viewed as the residue ring of the first by the ideal $(x_{2m-r},\ldots, x_{2m-1})$. 
We emphasize that the first is the ground polynomial ring of $\mathcal{H}_m$ and the second that of $\mathcal{H}_m[r]$.

\begin{Proposition}\label{subminors-are-cm}
	Assume that $0\leq r\leq m-2$ and let $1\leq t\leq m$. Then
	\begin{enumerate}
		\item[{\rm (i)}]  ${\rm codim}\, I_t(\mathcal{H}_m[r])=\min\{2(m-t)+1, 2m-t-r\}$ and  $\bar{R}/I_t(\mathcal{H}_m[r])$ is a Cohen--Macaulay ring.  
		\item[{\rm (ii)}]  $I_{t}(\mathcal{H}_m[r])$ is a prime ideal if and only if either $t=1$ or else $t\geq r+2$.
	\end{enumerate} 		
\end{Proposition}
\demo
(i)  We apply Lemma~\ref{homo1} (b), by which the ideal $I_{t}(\mathcal{H}_m[r])\subset \bar{R}$ is generated by the maximal minors of $\mathcal{H}_{t,2m-t}[r]$.
Therefore its codimension is at most $2m-t-t+1=2(m-t)+1$. 

We analyse the two cases separately:

{\sc (1) $ t \geq r+2$.}

Let $\mathcal{H}$ denote the uniquely defined generic Hankel matrix of size $t\times (2m-t)$  over the ground ring $R$ (see Remark~\ref{same_vars}). 
By Proposition~\ref{GP}, one has $I_t(\mathcal{H}_m)=I_t(\mathcal{H})$.
Setting $A:=R/I_t(\mathcal{H})$, one has $R/(x_{2m-r},\ldots, x_{2m-1})\simeq \bar{R}$ and 
$$A/(x_{2m-r},\ldots, x_{2m-1})A\simeq \bar{R}/I_t(\mathcal{H}_m[r]).$$
Since the strand $\{x_{2m-r},\ldots, x_{2m-1}\}$ is a subset of the system of parameters determined in Lemma~\ref{basic_minors}, as applied with $a=t, b=2m-t$, it is a regular sequence over $A$.
Therefore, $\bar{R}/I_t(\mathcal{H}_m[r])$ is Cohen--Macaulay.

Moreover, since $\dim(\bar{R}/I_t(\mathcal{H}_m[r]))=\dim(A)-r$, one gets
$\codim(I_{t}(H_m[r]))=2(m-t)+1$. 

\smallskip

{\sc (2) $ t \leq r+1$.}

In this situation, we get  $2m-t-(2m-r-1)=r-t+1$ null columns at the right end of the matrix $\mathcal{H}_{t,2m-t}[r]$, blurring the original requirement for the nature of the degeneration. Yet, removing these additional null columns yields a matrix of the shape $\mathcal{H}_{t,2m-r-1}[t-1]$ such that  $I_t(\mathcal{H}_{t,2m-r-1}[t-1])=I_{t}(\mathcal{H}_m[r])$. 

Applying \cite[Proposition 4.3]{Eisenbud2}  to the generic Hankel matrix $\mathcal H$ of size $t\times (2m-r-1)$ over the ring $R=k[x_1,\ldots, x_{2m-r+t-2}]$ yields that $R/I_t(\mathcal H)$ has dimension $2(t-1)$.
By the same token as before, the strand $\{x_{2m-r},\ldots, x_{2m-r+t-2}\}$ is a subset of the system of parameters determined in Lemma~\ref{basic_minors}, hence  it is a regular sequence over $A=R/I_t(\mathcal H)$.
Therefore, $\bar{R}/I_t(\mathcal{H}_m[r])$ is Cohen--Macaulay.
Doing the arithmetic once more, one gets ${\rm codim} \, I_{t}(H_m[r])=2m-t-r$.

\smallskip

(ii) Clearly, $I_1(\mathcal{H}_m[r])$ is prime.
Now, suppose that $t\leq r+1$. If $t=r+1$, the  $t$-minor on the bottom-right corner of the matrix contains a product of variables, hence is not prime. Clearly, this implies that no ideal of $t$-minors is prime either, for $t\leq r+1$.

For the reverse implication we proceed as follows.
Note that the hypothesis put us  in the first case of (i), in which the ideal $I_{t}(\mathcal{H}_m[r])$ coincides with the ideal of maximal minors of the matrix $\mathcal{H}_{t,2m-t}[r]$ and has codimension $2(m-t)+1$. 
Since the generic Hankel matrix is $1$-generic (see \cite[Proposition 4.3]{Eisenbud2}) then by \cite[Theorem 1 (ii)]{Eisenbud1}   the ideal $(x_{i_1},...,x_{i_s}, I_{t}(\mathcal{H}_{t,2m-t}))$ is prime provided $s\leq t-2$, where $\mathcal{H}_{t,2m-t}$ denotes the corresponding generic Hankel matrix.
But since $r\leq t-2$ by hypothesis, we are through.
\qed

\begin{Remark}\label{linear_span}\rm
	Note that a particular feature of the above result and its method of proof is that, for any given $1\leq t\leq m$ and any $0\leq r\leq m-2$, the dimension of the $k$-linear span of the $t$-minors of $\mathcal H_m[r]$ is the same as that of $\mathcal H_m$.
\end{Remark}

For convenience, we isolate two special cases of the above proposition.

\begin{Corollary}\label{degeneration_is_prime}
	Assume that $0\leq r\leq m-2$. One has:
	\begin{enumerate}
		\item[{\rm (i)}] If $m\geq 3$, $I_{m-1}(\mathcal{H}_m[r])$ has  codimension $3$ and is a prime ideal if and only if $r\leq m-3$.
		\item[{\rm (ii)}] If $m\geq 4$, $I_{m-2}(\mathcal{H}_m[r])$ has codimension $5$ and is a prime ideal if and only if $r\leq m-4$.
	\end{enumerate} 
\end{Corollary}

\begin{Remark}\rm
	It is interesting to note that, even for $r\leq m-3$, the ring $R/I_{m-1}(\mathcal{H}_m[r])$ is not always normal, a property that may require $m >> r$.
\end{Remark}

We add some additional considerations about the ideal of submaximal minors.
The next result is a non-generic version of \cite[Theorem 10.16 (b)]{B-V}, with the same proof.

\begin{Lemma}\label{subminors_spread}
	Let $\mathcal M$ be a square matrix whose entries are either variables over a field $k$ or zeros, such that $\det(\mathcal M)\neq 0$. Let $R$ denote the polynomial ring over $k$  on the nonzero entries of $\mathcal M$ and let $S\subset R$ denote the $k$-subalgebra generated by the submaximal minors.
	Then the extension $S\subset R$ is algebraic at the level of the respective fields of fractions.
\end{Lemma}

A consequence is the following:
\begin{Proposition}\label{max_anal_spread}
	With the notation of {\rm Proposition~\ref{subminors-are-cm}}, the ideal $I_{m-1}(\mathcal{H}_m[r])\subset R$ has maximal analytic spread.
\end{Proposition}
\demo The analytic spread of an ideal is also the dimension of the special fiber algebra of its Rees algebra. In this case, the ideal is a homogeneous ideal of the polynomial ring $R$ generated in one single degree. Thus, this algebra is isomorphic to the $k$-subalgebra $S\subset R$ generated by the minors.
Now apply the previous lemma.
\qed

\medskip

A tall order in classical invariant theory is the structure of the ideal of the defining polynomial relations of the submaximal minors of a matrix of linear entries. 
It may be interesting to compare the generic, generic symmetric and the Hankel situations as regards invariants of the coordinate section degeneration in the square case. Thus, let the size be $m\times m$ with $0\leq r\leq m-2$.
Let PI and CM be short for polar image and Cohen--Macaulay, respectively.
The overall picture looks as follows:

\medskip
\noindent
{\small
\begin{tabular}{|c|c|c|c|c|}
	\hline  & rank of Hessian & polar image (PI)  & primality of $I_{m-1}$ & fiber of    $I_{m-1}$\\ [2pt]
	\hline Generic & $m^2- r(r+1) $ & Gorenstein ladder  & if $m\geq {{r+1}\choose 2}+3 $  &  cone over PI\\ [2pt]
	\hline Symmetric & ${{m+1}\choose 2}-2\mathfrak{o}(r)$ & CM ladder &? & cone over PI\\ [2pt]
	\hline
	Hankel &$2m-r-1$ (maximum) & trivial  & if  $m\geq r+ 3$ & ?\\ [2pt]
	\hline  
\end{tabular}
}

\medskip

For the  first two rows of the above table see \cite{Degen_generic} and \cite{Degen_sym}, respectively.
In regard to the question mark at the end of the third row, one can reduce to the case where $I_{m-1}$ is replaced by the maximal minors of an $(m-1)\times (m+1)$ degenerate Hankel matrix as stated in Lemma~\ref{homo1} (b) and dealt with in the previous subsection. In the case where $r=0$ (generic Hankel), one knows that the defining ideal of the fiber is generated by Pl\"ucker relations (\cite[Theorem 4.7]{Conca1998}) -- more precisely, the fiber is isomorphic as $k$-algebra to the homogeneous coordinate ring of the Grassmannian $G(m-1,m+1)$ as is proved in Proposition~\ref{Hankel_is_Grassmann}.
However, already for $m=4, r=1$, there is a minimal defining relation of degree $3$ besides the quadratic Pl\"ucker relations.

This state of affairs may lead to the following:

\begin{Question}\rm Let $1\leq r\leq m-2$.
	\begin{itemize}
		\item Is the special fiber of  $I_{m-1}(\mathcal{H}_m[r])$ Cohen--Macaulay?
		\item Is the Rees algebra of $I_{m-1}(\mathcal{H}_m[r])$ Cohen--Macaulay and of fiber type?
		\item Is the defining ideal of the special fiber of $I_{m-1}(\mathcal{H}_m[r])$  minimally generated by quadrics and cubics? 
		\item Are the cubic ones merely forced by the degeneration assumption or have a deeper $GL$-representation meaning as described in \cite{BCV}?
	\end{itemize}
\end{Question}

These  questions are specially intriguing because they actually refer to maximal minors.
For low values of $m$, a computer verification has been carried for the first three questions.
Curiously, in the case $r=m-2$ the defining ideal of the special fiber has even minimal binomial generators with  non squarefree terms.

\section{The Hessian and the Gradient ideal}\label{Hessian}
We emphasize that throughout  the ground field $k$ has characteristic zero (or sufficently large as compared to the size $m$ of the matrix).

The gradient ideal $J=J(f)\subset R$  of $f:=\det \mathcal{H}_m[r]$ is the ideal generated by the partial derivatives of $f$.
The Jacobian matrix of the partial derivatives is the so-called {\em Hessian matrix} of $f$ and its determinant is called {\em the Hessian} of $f$.

The basic result is

\begin{Theorem}\label{Hessian_not_zero} Let $f=\det \mathcal H{_m}[r]$.
For $0\leq r\leq m-2$, the Hessian $h(f)$ does not vanish.
\end{Theorem}
Due to its length, the proof is given in the Appendix.
An argument  in the spirit of  Corollary~\ref{Hessian_is_not_null}, obtained in the  fully generic case, is more than welcome, but so far we have been unable to hit it.

\subsection{Parabolism: the determinant versus its Hessian}

It is a classical question as to when a form $f$ is a factor of its Hessian $h(f)$ with multiplicity $\geq 1$. 
If this takes place, then $f$ is said in addition to have the {\em expected multiplicity} (according to Segre) if its multiplicity as a factor of $h(f)$ is  $ {\rm codim}(V(f)^*)-1$, where $V(f)^*$ denotes the dual variety to the hypersurface $V(f)$ (see \cite{CRS}).

In this part we elaborate on the case of $f:=\det {\mathcal H}_m[r]$. 
 The following is a first step regarding the Segre expected multiplicity property.

\begin{Proposition}\label{dim_dual_lower_bound}
	 Let $0\leq r\leq m-2 \, (m\geq 3),$ and $f=\det \mathcal H_m[r]$. 
		Then $\dim V(f)^* $ is at least $m-1$.	
\end{Proposition}
\demo 
A formula due to Segre, as transcribed in \cite[Lemma 7.2.7]{Frusso}, says that 
$$\dim (V(f)^*)={\rm rank }\; H(f)({\rm mod}\, f)-2,$$
where $H(f)$ denotes de Hessian matrix $f$. It will then suffice to show that $H(f)$ has rank at least $m+1$ modulo $f$ (note, as a slight control, that $m+1\leq 2m-r-1$ for $r\leq m-2$). 

Consider the Hankel matrix $\mathcal{H}_m(m-2)$ obtained by further degenerating 
$$x_{m+2}= \cdots = x_{2m-(r+1)}=0.$$
The Hessian matrix of $\det \mathcal{H}_m(m-2)$  can be viewed as the $(m+1)\times (m+1)$ submatrix $\Theta$  of $H(f)$ of the first $m+1$ rows and columns modulo   $(x_{m+2}, \ldots , x_{2m-(r+1)} )$.
By \cite[Theorem 4.4(iii)]{CRS}, $\det \Theta$ modulo   $(x_{m+2}, \ldots , x_{2m-(r+1)} )$ is a nonzero scalar multiple of  $x_{m+1}^{(m+1)(m-2)}$. Thus, $\det \Theta$ does not vanish. In addition, $\det \Theta$ does not vanish module $f$ because $f$ does not have a pure term in $x_{m+1}.$ Therefore, $H(f)$ has rank at least $m+1$ modulo $f$.
\qed

\smallskip

The context immediately prompt us to a few essential questions, which we choose to state as:

\begin{Conjecture}\label{dual_variety} Let $0\leq r\leq m-2 \, (m\geq 3),$ and $f=\det \mathcal H_m[r]$. 
	Then:
	\begin{enumerate}
		\item[{\rm (i)}] $\dim V(f)^*  = m-1$, and hence the expected multiplicity of $f$ as a factor of its Hessian determinant is $2m-r-2-(m-1)-1=m-r-2$.	
		\item[{\rm (ii)}] $f$ is a factor of its Hessian determinant with the expected multiplicity.	
		\item[{\rm (iii)}]  Let $\mathcal D$ denote the defining ideal of $V(f)^*$ in its natural embedding. Then the initial degree of $\mathcal D$ is $m$ and the subideal $(\mathcal D_m)$ generated in the initial degree has codimension $m-r-1$.
		{\rm (}In addition, if $0\neq r=m-3$ then $(\mathcal D_m)$ is a  linearly presented codimension $2$ perfect ideal.{\rm )}
		\item[{\rm (iv)}] $V(f)^*$ is arithmetically Cohen--Macaulay if and only if  $r=m-2$, in which case $V(f)$ is self-dual up to a coordinate change.
	\end{enumerate}
\end{Conjecture} 

Since a conjecture ought to be grounded on reasons other than mere computer experimentation, we give some elements for potential proofs:

$\bullet$ The conjectured statement in (i) accommodates the generic case as well as the sub-Hankel degeneration ($r=m-2$) -- for the latter,  $f$ is not a factor of its Hessian determinant (cf. \cite[Theorem 4.4 (iii)]{CRS}).
	
$\bullet$ The value of the dimension in item (i) follows from  Proposition~\ref{dim_dual_lower_bound} and the conjectured codimension in item (iii).
Note that (iii) tells, in particular, that the Hankel degenerations of the sort considered here have geometric behavior totally distinct from the fully generic and the symmetric counterparts -- for the latter the dual variety is defined by quadrics and often by ladder $2$-minors (see \cite{Degen_generic} and \cite{Degen_sym}).

$\bullet$ As for (iv), one has an affirmative answer in one direction, as follows. Thus, suppose that $r=m-2$ and let us show that $V(f)^*$ is self-dual up to coefficients (in particular, $V(f)^*$ is arithmetically Cohen--Macaulay).

Set $\bar{R}:=R/f$ and $\bar{J}:=(J,f)/f$ where $J$ is the gradient ideal of $f.$ By \cite[Lemma 4.2]{CRS} and Euler's formula   the syzygy matrix of $\bar{J}$ contains a $(m+1)\times(m+1)$ submatrix with the following shape:

$$\varphi:=\left(\begin{array}{cccccccccccccc}
x_1&(m-1)x_1&2x_2&2x_3&2x_4&\ldots&2x_{m}\\
x_2&(m-2)x_2& \ast x_3&\ast x_4&\ast x_5&\ldots&x_{m+1}\\
x_3&(m-3)x_3&\ast x_4&\ast x_5&\ast x_6&\ldots &0\\
\vdots&\vdots&\vdots&\vdots&\vdots&\cdots&\vdots\\
x_{m-3}&3x_{m-3}&\ast x_{m-2}&\ast x_{m-1}&\ast x_m&\ldots&0\\
x_{m-2}&2x_{m-2}&\ast x_{m-1}&\ast x_m&x_{m+1}&\ldots&0\\
x_{m-1}&x_{m-1}&\ast x_m&x_{m+1}&0&\ldots&0\\
x_m&0&x_{m+1}&0&0&\ldots&0\\
x_{m+1}&-x_{m+1}&0&0&0&\ldots&0
\end{array}\right) \,(\bmod f)$$
where $\ast$ are unspecified nonzero coefficients. Since $\phi$ is a submatrix of the syzygy matrix, its rank is at most $m.$ Expanding the determinant of $\phi$ by Laplace along the last row, one gets
$x_{m+1}\det \widetilde{\phi}\equiv 0 \,(\bmod f)$, where
$$\widetilde{\phi}=\left(\begin{array}{cccccccccccccc}
mx_1&2x_2&2x_3&2x_4&\ldots&2x_{m}\\
(m-1)x_2& \ast x_3&\ast x_4&\ast x_5&\ldots&x_{m+1}\\
\vdots&\vdots&\vdots&\vdots&\cdots&\vdots\\
4x_{m-3}&\ast x_{m-2}&\ast x_{m-1}&\ast x_m&\ldots&0\\
3x_{m-2}&\ast x_{m-1}&\ast x_m&x_{m+1}&\ldots&0\\
2x_{m-1}&\ast x_m&x_{m+1}&0&\ldots&0\\
x_m&x_{m+1}&0&0&\ldots&0\\
\end{array}\right).$$
Since $x_{m+1}\not\equiv  0 \,(\bmod f)$ and  $\bar{R}$ is a domain we have $\det\widetilde{\phi}\in(f).$ Note that $\det \widetilde{\phi}\neq 0.$  In particular, for reasons of degree,   $\det \widetilde{\phi}$ is a nonzero scalar multiple of  $f.$

Let $B$ be the unique  $(m+1)\times(m+1)$  linear matrix with entries in $k[t_1,\ldots,t_{m+1}]=k[\tt]$ such that $$(\tt)\cdot\varphi=(\xx)\cdot B.$$ In particular, $B$ has the following shape:
$$\left(\begin{array}{ccccccccccccccccccc}
t_1&(m-1)t_1&0&0&0&\ldots&0\\
t_2&(m-2)t_2&2t_1&0&0&\ldots&0\\
t_3&(m-3)t_3&\bullet t_2&2t_1&0&\ldots&0\\
t_4&(m-4)t_4&\bullet t_3&\bullet t_2&2t_1&\ldots&0\\
\vdots&\vdots&\vdots&\vdots&\vdots&\cdots&\vdots\\
t_m&0&\bullet t_{m-1}&\bullet t_{m-2}&\bullet t_{m-3}&\ldots&2t_1\\
t_{m+1}&-t_{m+1}&t_m&t_{m-1}&t_{m-2}&\ldots&t_2
\end{array}\right)$$
where $\bullet$ are unspecified nonzero coefficients. Note that $B^t$ is  submatrix of the matriz Jacobian dual of $\bar{J}.$

Let $k[t_1,\ldots,t_{m+1}]/P$ denote the homogeneous coordinate ring of $V(f)^*$. By \cite{CRS} one knows that the Hessian matrix $H(f)$ has rank $m+1$ modulo $(f)$, hence $P$ is a principal (prime) ideal.

 Since the rank of the Jacobian dual matrix of $\bar{J}$ modulo $P$ is at most $m$ we have $\det B\in P.$ 
 Expanding $\det B$ by Laplace along the first row yields  $\det B=t_1\det\widetilde{B}$, where 
$$\widetilde{B}=\left(\begin{array}{ccccccccccccccccccc}
-t_2&2t_1&0&0&\ldots&0\\
-2t_3&\bullet t_2&2t_1&0&\ldots&0\\
-3t_4&\bullet t_3&\bullet t_2&2t_1&\ldots&0\\
\vdots&\vdots&\vdots&\vdots&\cdots&\vdots\\
-(m-1)t_{m}&\bullet t_{m-1}&\bullet t_{m-2}&\bullet t_{m-3}&\ldots&2t_1\\
-mt_{m+1}&t_m&t_{m-1}&t_{m-2}&\ldots&t_2
\end{array}\right).$$
In particular, $t_1\det\widetilde{B}\equiv 0\, (\bmod\, P).$ 
Since $P$ is prime and $t_1\notin P$, then $\det\widetilde{B}\in P.$  
On the other hand, the same sort of argument as \cite[Remark 4.6 (c)]{CRS} shows that  $ \det \widetilde{B} $ is irreducible. Thus, $P=(\det\widetilde{B}).$
To conclude, note that $\widetilde{B}$ can be replaced by the obvious sub-Hankel matrix with nonzero coefficients obtained by exchanging rows equidistant from the extremes.

\subsection{The gradient ideal}

\subsubsection{The codimension}

The codimension of the gradient ideal $J$ of $\det(\mathcal{H}_m[r])$ comes in as one of the basic ingredients in the paper.
The result below is a neat consequence of the nature of $J$ in its relation to the cofactors of $\mathcal{H}_m[r]$ and of the elementary cofactor formulas.

\begin{Proposition}\label{cod_of_gradient}{\rm (char$(k)=0$)}
Let $J\subset R=k[x_1,\ldots,x_{2m-r-1}]$ denote the gradient ideal of $\det(\mathcal{H}_m[r])$, where $m-r\geq 2$.
Then
$$ {\rm codim}(J)= \left\{
\begin{array}{ll}
2 & \mbox{ if $m-r=2$}\\
3 & \mbox{otherwise.}
\end{array}
\right.
$$ 
\end{Proposition}
\demo 
By Proposition~\ref{partials_vs_cofactors}, $J$ is contained in $I_{m-1}(\mathcal{H}_m[r])$, for every degeneration step. The latter has codimension $3$ by Corollary~\ref{degeneration_is_prime}.
Therefore, $J$ has codimension  at most $3$.

The case where  $m-r=2$ is easily checked and, in any case, sufficiently studied in \cite{CRS} and \cite{MAron}.
Thus, we assume that $m-r\geq 3$ and induct on $m-r$.

When $m-r=3$, one proceeds as follows. 

Let $\Delta_{i,j}$ denote the (signed) cofactor of the $(j,i)$-entry of $\mathcal{H}_m[r].$ Given a prime ideal $P$ containing $J$, we will prove that $P\supset I_{m-1}(\mathcal{H}_m[r])$ or $P\supset (x_m,x_{m+1},x_{m+2})$, which implies that any minimal prime of $J$ has codimension at least $3$. For this, we divide the proof in two cases.

\medskip
{\sc Case 1:} $x_{m+2}\in P.$

Note that $f_1= x_{m+1}^{m-1}+H$ with $H\in (x_{m+2})$. Thus, since $(x_{m+2}, f_1)\subset P $ then $x_{m+1}\in P.$ By a similar token, $f_m=x_{m}^{m-1}+G$ with $G\in (x_{m+1},x_{m+2}).$ This way, since $(x_{m+1},x_{m+2}, f_{m})\in P,$ we see that $x_{m}\in P.$ Therefore, $(x_m,x_{m+1},x_{m+2})\subset P.$

\medskip

{\sc Case 2:} $x_{m+2}\notin P.$

We claim that, for every $2\leq k\leq m$, the entries of the matrix

$$\left( \begin{array}{cccc}
\Delta_{1,1} & \cdots & \Delta_{1,k} \\
\vdots & \cdots & \vdots	\\
\Delta_{k,1} &  \cdots& \Delta_{k,k} 
\end{array}\right)$$
belong to $P.$ 

Induct on $k$.

Let $k=2$.
As a consequence of the classical formula for the matrix of cofactors of the above matrix one has the following relation:
\begin{equation}\label{kequal2}
x_{m}\Delta_{1,1}+x_{m+1}\Delta_{2,1}+x_{m+2}\Delta_{3,1}=0
\end{equation}
Since $f_1=\Delta_{1,1}$ and $f_2=(1/2)\Delta_{2,1},$ it follows from \eqref{kequal2} that $\Delta_{3,1}\in P.$ Thus, since  $f_3=\Delta_{2,2}+2\Delta_{3,1}$ we have $\Delta_{2,2}\in P.$ Hence, for $k=2$ the statement is true.

Now, assuming the statement for some let $k\geq  2$, we show it holds for $k+1.$ Consider the matrix

$$\left( \begin{array}{ccc|c}
\Delta_{1,1} & \cdots & \Delta_{1,k} & \Delta_{1,k+1}\\
\vdots & \cdots & \vdots & \vdots	\\
\Delta_{k,1} &  \cdots& \Delta_{k,k} & \Delta_{k,k+1}\\
\hline
\Delta_{k+1,1} &  \cdots& \Delta_{k+1,k} & \Delta_{k+1,k+1}
\end{array}\right) $$
Again  the cofactor formula yields the following relations:
\begin{equation}\label{eq-cod}
x_{m-k+2}\Delta_{1,j}+\cdots +x_{m+1}\Delta_{k,j}+x_{m+2}\Delta_{k+1,j}=0,\  \mbox{for all } 1\leq j\leq k+1.
\end{equation}
These equalities, along with the inductive hypothesis, yield 
\begin{equation}\label{last_Delta}
\Delta_{j,k+1},\Delta_{k+1,j}\in P
\end{equation}
for every $1\leq j\leq k.$ 
Finally, for $j=k+1$ one has 
\begin{equation*}
x_{m-k+2}\Delta_{1,k+1}+\cdots +x_{m+1}\Delta_{k,k+1}+x_{m+2}\Delta_{k+1,k+1}=0
\end{equation*}
From this and \eqref{last_Delta}, it follows that $\Delta_{k+1,k+1}\in P.$ 

Therefore, taking $k=m$ we have $I_{m-1}(\mathcal{H}_m[r])\subset P$, as was to be shown.

Thus, we are through with the case where $m-r=3$. 

We now induct on $m-r$.
For the inductive step, note that the ascending induction step from $m-r$ to $m-r+1=m-(r-1)$ corresponds to a descending induction step from $r$ to $r-1$. 
Thus, we are given the matrix $\mathcal{H}_m[r-1]$, with $r-1\leq m-4$, and the corresponding gradient ideal $J\subset R=k[x_1,\ldots,x_{2m-r}]$ and assume by induction that the corresponding gradient ideal $J'\subset R'=k[x_1,\ldots,x_{2m-r-1}]$ of $\det (\mathcal{H}_m[r])$ has codimension $3$.
Since the latter matrix is a degeneration of the former by setting $x_{2m-r}\mapsto 0$, clearly $(J',x_{2m-r})\subset (J,x_{2m-r})$ as ideals in the bigger ring $R$, hence the codimension of $J$ is at least $3$ since $(J',x_{2m-r})$ has codimension $4$. 
\qed

\begin{Remark}\label{odds4conjecture}\rm
	(1) The proposition will be a consequence of the more encompassing Theorem~\ref{prime_ideal_structure} (a) below. 
	Indeed, the proofs bear some similarity, but are mutually independent.
	 
	(2)
It has been seen in the proof of Proposition~\ref{subminors-are-cm} (i) that, for arbitrary $r\leq m-3$, the entry subset $\{x_{2m-r}, x_{2m-(r-1)},\ldots, x_{2m-1}\}$ of the fully generic $m\times m$ Hankel matrix $\mathcal{H}_m$ is a regular sequence modulo the submaximal minors of $\mathcal{H}_m$. {\em It would seem} that this sequence is equally regular modulo the gradient ideal of $\det(\mathcal{H}_m)$ (see Conjecture~\ref{exact_sequence} below). Thus,  the specialization still has codimension $3$ in the entry ring of $\mathcal{H}_m[r]$. Unfortunately, the specialized ideal is larger than the gradient ideal of $\det(\mathcal{H}_m[r])$, hence one resorts to the inductive argument above to bypass this apparent obstruction.
\end{Remark}

\subsubsection{The associated primes}

In this part we will suppose throughout that $r\leq m-3$.
The case where $r=m-2$ has been thoroughly dissected in \cite{CRS} and \cite{MAron}.

We first develop a few ideas around lower ideals of minors and the cofactor formula, some of which have been scratched in the previous section. 

Fix an integer $1\leq j\leq m-2.$

Consider block partitions of $\mathcal{H}_m[r]$ and its cofactor matrix, as follows: 
\begin{equation}\label{decompose_Hankel}
\mathcal{H}_m[r]=\left[\begin{array}{cc}
U_{m-j}\\
D_j
\end{array}\right]
\end{equation}
where $U_{m-j}$ is an $m-j$-rowed rows.

\begin{equation}\label{decompose_adj}
{\rm cof}(\mathcal{H}_m[r])=\left[\begin{array}{cc}
A_{m-j}&B_{m-j}\\
B'_j&C_j
\end{array}\right], \quad B'_j=B_{m-j}^t
\end{equation}
where $A_{m-j}$ is a square $(m-j)$-rowed matrix.

The notation is such that the subscript of any of the blocks denotes its number of rows.

Block multiplication and the cofactor formula yield
\begin{equation}\label{block_cofactor}
f\, \mathbb{I}_m={\rm cof}(\mathcal{H}_m[r])\, \mathcal{H}_m[r]=
\left[\begin{matrix}
A_{m-j} U_{m-j}+B_{m-j} D_j\\
B_{m-j}^t U_{m-j}+C_j D_j
\end{matrix}
\right],
\end{equation}
where  $f=\det \mathcal{H}_m[r].$
Since char$(k)=0$, Euler's identity implies that $f\in J=J_m[r]\subset R=k[x_1,\ldots,x_{2m-r-1}]$, and hence the entries on the rightmost matrix belong to $J$ as well.

\begin{Lemma}\label{contemprod}
	For any ideal $I\subset R$ containing $J$, the following holds:
	$$I_1(A_{m-j}U_{m-j})\subset I\Rightarrow I_1(B_{m-j})I_j(D_j)\subset I.$$
	By the same token, 
	$$I_1(B_{m-j}^t U_{m-j})\subset I\Rightarrow I_1(C_j)I_j(D_j)\subset I.$$	
%	Se $I$ \'e um ideal que cont\'em $J$ e as entradas de $A_{m-j}\cdot U_{m-j}$ $($resp., $B_{m-j}^t\cdot U_{m-j}$$)$ ent\~ao $I$ cont\'em $I_{j}(D_j)I_1(B_{m-j})$ $($resp., $I_{j}(D_j)I_1(C_j)$$)$
\end{Lemma}
\demo It suffices to consider the first implication.
As has been noted above, one has $I_1(A_{m-j} U_{m-j}+B_{m-j} D_j)\subset J$, hence  $I_1(A_{m-j} U_{m-j}+B_{m-j} D_j)\subset I$ and since $I_1(A_{m-j}U_{m-j})\subset I$ by hypothesis, then $I_1(B_{m-j}D_j)\subset I$.

Note that $I_j(D_j)$ is the ideal of maximal minors of $D_j$.
Thus, letting $M$ denote an arbitrary $j\times j$ submatrix of $D_j$, one has
 $I_1(B_{m-j} M)\subset I$. Thus,  $I_1(B_{m-j}\, M\, {\rm cof}(M))= I_1(\det(M)\, B_{m-j})\subset I$. 
 Therefore, $I_{j}(D_j)I_1(B_{m-j})\subset I$,
as required. 
\qed

\smallskip

Note that Lemma~\ref{contemprod}  still holds by replacing the four block matrices of cof$(\mathcal{H}_m[r])$ by their respective $i$th rows, for an arbitrary $i$.

Given a matrix $\mathcal{M}$ and an index $i$, $L_i(\mathcal{M})$ will stand for its $i$th row.

For the next lemma, write $[L_1(A_{m-j})\,\, L_1(B_{m-j})]$  for the first row of the matrix ${\rm cof}(\mathcal{H}_m[r])$ as decomposed in (\ref{decompose_adj}).
Recall a previous convention by which $\Delta_{i,j}$ denotes the (signed) cofactor of the $(j,i)$-entry of $\mathcal{H}_m[r].$

\begin{Lemma}\label{lema_mmenos2} {\rm ($j=m-2$)} $I_{m-2}(D_{m-2})\cdot I_1([L_1(A_{2})\,\, L_1(B_{2})])\subset J.$
\end{Lemma}
\demo Clearly,  $L_1(A_{2})=[\Delta_{11}\,\,\Delta_{1,2}],$ hence its entries belong to $J$ since  $\partial f/\partial x_1=\Delta_{1,1}$ and $\partial f/\partial x_2=2\Delta_{1,2}.$ Clearly, then $I_1(L_1(A_{2}) U_{2})\subset J$. By Lemma~\ref{contemprod} one deduces that
$I_{m-2}(D_{m-2})\cdot I_1(L_1(B_{2}))\subset J.$
\qed

\smallskip

So far, $j$ was fixed. Now we make it vary in a certain range in order to decide when a given prime ideal containing $J$ contains or not the ideal $I_{m-1}(\mathcal{H}_m[r])$.
For example, the prime ideal $Q=(x_m,x_{m+1},\ldots,x_{2m-r-1})$ contains $J$ but not $I_{m-1}(\mathcal{H}_m[r])$, while adding the entry $x_{m-1}$ to $Q$ gives a prime ideal containing the submaximal minors.
Precisely, one has:

\begin{Proposition}\label{lema_ind} Suppose that $r\leq m-4$.
	Let  $P\supset J$ be a prime ideal. If there exists an index $j$ in the range $r+2\leq j\leq m-2$ such that  $I_{j}(D_j)\subset P$ and $I_{j-1}(D_{j-1})\not\subset P$, then $I_{m-1}(\mathcal{H}_m[r])\subset P.$
\end{Proposition}
\demo Since $I_{j}(D_j)\subset P$, (\ref{block_cofactor}) implies that the entries on the row $[A_{m-j}\,\,B_{m-j}]$ (resp., $[A_{m-j}\,\,B_{m-j}]^t$) belong to $P.$  On the other hand, one has:

$$\Delta_{m-j+1,m-j+1}=\partial f/\partial x_{2(m-j+1)-1}-2\sum_{i=\max\{1,m-2j+2\}}^{m-j} \Delta_{2(m-j+1)-i,i} \in P.$$

Therefore, the row entries of $L_{m-j+1}(A_{m-j+1})=[\Delta_{m-j+1,1}\ldots\Delta_{m-j+1,m-j+1}]$ belong to $P.$ It follows that $I_1(L_{m-j+1}(A_{m-j+1})\cdot U_{m-j+1})\subset P.$ 
Then the upshot from Lemma~\ref{contemprod} is that  
$$I_{j-1}(D_{j-1})\cdot I_1(L_{m-j+1}(B_{m-j+1}))\subset P.$$ 
But since  $I_{j-1}(D_{j-1})\not\subset P,$ then  $$I_1(L_{m-j+1}(B_{m-j+1}))=(\Delta_{m-j+1,m-j+2},\ldots,\Delta_{m-j+1,m})\subset P.$$
This in turn tells us that  

\begin{equation}\label{AB}
I_1([A_{m-j+1}\,\,B_{m-j+1}])\subset P \quad\quad (\mbox{resp.},\,\, I_1([A_{m-j+1}\,\,B_{m-j+1}]^t)\subset P)
\end{equation}
Thus, $I_1(B_{m-j+1}^t\cdot U_{m-j+1})\subset P$ and once again Lemma~\ref{contemprod}, implies that  
$$I_{j-1}(D_{j-1})I_1(C_{j-1})\subset P.$$
It follows that
\begin{equation}\label{C}
I_1(C_{j-1})\subset P.
\end{equation}
A moment reflection will convince us that \eqref{AB} and \eqref{C} imply that $I_{m-1}(\mathcal{H}_m[r])\subset P$, as was to be shown.
\qed

\begin{Theorem}\label{prime_ideal_structure}
Let $J\subset R$ denote the gradient ideal of the determinant of $\mathcal{H}_m[r]$, with $1\leq r\leq m-3$, and let $Q$ denote the ideal generated by the $m-r$ nonzero variables of its last column. 
Then:
\begin{enumerate}
	\item[{\rm (a)}] The minimal primes of $R/J$ are $Q$ and $P:=I_{m-1}(\mathcal{H}_m[r])$.
	\item[{\rm (b)}] $J$ is  not a reduction of $P$.
	\item[{\rm (c)}] The unmixed and minimal components of $J$ coincide if and only if $r=m-3$.
\end{enumerate}
\end{Theorem}
\demo 
 (a) It suffices to prove that any prime ideal $P$ containing $J$ necessarily contains either $Q=(x_m,x_{m+1},\ldots,x_{2m-r-1})$ or $I_{m-1}(\mathcal{H}_m[r]).$ 
 Divide the argument in two cases:

 \medskip
 
 \noindent {\bf Case 1:} $I_{r+1}(D_{r+1})\subset P.$
 
 Since $x_{2m-r-1}^{r+1}\in I_{r+1}(D_{r+1})$ then $x_{2m-r-1}\in P.$  More: for any $m\leq u\leq 2m-r-2$ there exists an $(r+1)$-minor in $I_{r+1}(D_{r+1})$ of the form
 $$x_u^{r+1}+x_{u+1}g_{u+1}+\cdots+x_{2m-r-1}g_{2m-r-1},$$ 
 for certain forms $g_{u+1},\ldots,g_{2m-r-1}\in R$.
 Decreasing induction on $u$ then wraps up the inclusion $(x_m,\ldots,x_{2m-r-1})\subset P.$ 
 
 \medskip
 
 \noindent {\bf Case 2:} $I_{r+1}(D_{r+1})\not\subset P.$
 
 If there exists an index $j$ in the range $r+2\leq j\leq m-2$ such that $I_j(D_j)\subset P$ then Proposition~\ref{lema_ind} implies that  $I_{m-1}(\mathcal{H}_m[r])\subset P.$ 
 
 Thus, assume that $I_j(D_j)\not\subset P$ for every $r+2\leq j\leq m-2.$ In particular, $I_{m-2}(D_{m-2})\not\subset P.$ By Lemma~\ref{lema_mmenos2},
 $$I_{m-2}(D_{m-2})\cdot I_1([L_1(A_{2})\,\, L_1(B_{2})])\subset P$$
 and since $I_{m-2}(D_{m-2})\not\subset P,$ then 
 \begin{equation}\label{inclusao1}
 I_1([L_1(A_{2})\,\, L_1(B_{2})])\subset P.
 \end{equation} 
 Thus, $\Delta_{1,3}\in P.$ But $\Delta_{2,2}=\partial f/\partial x_3-2\Delta_{1,3}\in J\subset P.$ Therefore, the entries  $L_2(A_{2})=[\Delta_{2,1}\,\,\Delta_{2,2}]$ belong to $P.$ 
 Then, from Lemma~\ref{contemprod} one has 
 $$I_{m-2}(D_{m-2})\cdot I_1([L_2(A_{2})\,\, L_2(B_{2})])\subset P$$ 
 hence, 
 \begin{equation}\label{inclusao2}
 I_1([L_2(A_{2})\,\, L_2(B_{2})])\subset P.
 \end{equation}
 From \eqref{inclusao1} and  \eqref{inclusao2} one has 
 \begin{equation}\label{AB2}
 I_1([A_{2}\,\,B_{2}])\subset P \quad\quad (\mbox{resp.},\,\, I_1([A_{2}\,\,B_{2}]^t)\subset P)
 \end{equation} 
 
 From these inclusions and  Lemma~\ref{contemprod}  it obtains $I_{m-2}(D_{m-2})\cdot I_{1}(C_{m-2})\subset P,$ and hence, 
 \begin{equation}\label{inclusao3}
 I_{1}(C_{m-2})\subset P.
 \end{equation}
 Clearly, \eqref{AB2} and \eqref{inclusao3} imply that $I_{m-1}(\mathcal{H}_m[r])\subset P.$ 
 
 \medskip

(b) If $J$ is a reduction of $P$, at least $\sqrt{J}=P$, which would contracdict the result in (a).

\medskip

(c) This is an immediate consequence of (a).
\qed

\begin{Remark}\label{Associated_primes}\rm
(i) Computational evidence suggests that the $Q$-primary component of $J$ is generated by ${{m-1}\choose r}$ forms of degree $r$ (but only coincides with the $r$th power of $Q$ when $r=1$).	

(ii) The structure of the embedded associated primes of $R/J$ is quite involved. The following two primes seem to be candidates in general: $(x_{m-1},Q)$ {\rm (}in codimension $m-r+1${\rm )} and $\sqrt{I_{m-2}(\mathcal{H}_m[r]})$ {\rm (}in codimension $5${\rm )}.
\end{Remark}

\subsubsection{Linear behavior}

We once more emphasize that the ground characteristic is assumed to be zero.

For the generic Hankel matrix $\mathcal{H}_m$ it was proved that the linear rank of $J$ is $3$ (\cite[Theorem 3.3.5]{Maral}).
At the other end of the spectrum, for the degeneration $\mathcal{H}_m[m-2]$ it was proved that the linear rank of $J$ is maximal possible ($=m$) (\cite{CRS}) and $J$ is in addition of linear type (\cite[Section 4.1, Theorem 4.8]{MAron}).

Though not obvious at all, one expects that the linear rank of $J_m[r]$ for $1\leq r \leq m-3$ be squeezed in-between.
 In this regard, one has two strong conjectures:
 
 \begin{Conjecture}\label{exact_sequence}
 If $J=J_m[0]$ is the gradient ideal of the fully generic $m\times m$ Hankel matrix, then the sequence $\{x_{m+3},\ldots, x_{2m-1}\}$ is regular modulo $J$.
 \end{Conjecture}

\begin{Conjecture}\label{linear_rank_is_2}
Let $1\leq r\leq m-3$. Then the linear rank of $J_m[r]$ is $2$.
\end{Conjecture}

{\sc Claim.} Conjecture~\ref{exact_sequence} implies Conjecture~\ref{linear_rank_is_2}.

\indent \demo  Set $R:=k[x_1,\ldots, x_{2m-1}]$,
 $\bar{R}:=R/(x_{2m-r},\ldots, x_{2m-1}),\; \bar{J}:=J\bar{R}\subset \bar{R},$
so
$$\bar{J}=(J,x_{2m-r},\ldots, x_{2m-1})/(x_{2m-r},\ldots, x_{2m-1})\subset k[x_1,\ldots,x_{2m-r-1}].
$$
One asserts that the module of linear syzygies of $\bar{J}$ is a free $\bar{R}$-module of rank $3$.

To see this, note that since $r\leq m-3$ then $2m-r\geq m+3$,  Conjecture~\ref{exact_sequence} implies that $\{x_{2m-r},\ldots, x_{2m-1}\}$ is a regular sequence modulo $J$.

The minimal graded resolution of  $R/J$ specializes to that of $\bar{R}/\bar{J}$. In particular, the module of linear syzygies of $\bar{J}$  has the same structure as that of $J$.
But, in the fully generic case, this module has been shown to be free of rank $3$ (\cite{Maral}).

Now, the generators of $J_m[r]$ are part of a minimal set of generators of $\bar{J}$ in the natural order of the partial derivatives, as follows from Proposition~\ref{partials_vs_cofactors}. In particular, any linear syzygy of $J_m[r]$ gives one of $\bar{J}$ by filling down enough zeros.
This implies that the module of syzygies of  $J_m[r]$ is a submodule of that of $\bar{J}$ and is free of rank $\leq 3$.

On the other hand, according to \cite[Chapter 3, Section 3.3.2, p.22]{Maral}, the submodule of linear syzygies of the fully generic $J$ has the matrix form
\begin{equation}\label{generic_syzygies}
\left(
\begin{matrix}
0 & \beta_1x_2 & \gamma_1x_1\\
\alpha_2x_1 &\beta_2x_3 &\gamma_2x_2\\
\vdots & \vdots & \vdots\\
\alpha_{2m-2}x_{2m-3} & \beta_{2m-1}x_{2m-1} &\gamma_{2m-2} x_{2m-2}\\
\alpha_{2m-1}x_{2m-2} & 0 & \gamma_{2m-1}x_{2m-1}
\end{matrix}
\right)
\end{equation}
with $\alpha_i,\beta_j,\gamma_l$ elements of $k$, where $\alpha_{i}\neq 0$, for $2\leq i\leq 2m-1$.
Note that there is no $k$-elementary operation that kills the last coordinate of the first syzygy above.
Since the module specializes to the module of linear syzygies of $\bar{J}$, the latter is obtained by setting to zero the entries $x_{2m-r},\ldots, x_{2m-1}$.
Thus, there is a linear syzygy of $\bar{J}$ whose $(2m-r)$th coordinate $\alpha_{2m-r}x_{2m-r-1}$ does not vanish.
Since $f_m[r]$ has only $2m-r-1$ derivatives, this syzygy cannot be a syzygy of $J_m[r]$.
Therefore, the module of syzygies of $J_m[r]$ is free of rank at most $2$.

But, by the same token, the $\beta$ and $\gamma$ version of the syzygies of $\bar{J}$ are syzygies of $J_m[r]$.
We conclude that the module of syzygies of $J_m[r]$ is free of rank exactly $2$; in particular, it has linear rank $2$.
\qed

\begin{Remark}\rm
	(1) Conjecture~\ref{linear_rank_is_2} fails in positive characteristic. For example, in characteristic $3$ the linear rank of $J_4[1]$ is $3$. The reason might be the juggling with the nonvanishing syzygies coordinates of (\ref{generic_syzygies}) made possible in null characteristic. 
	
(2) If $r=m-2$ a similar argument crumbles down because the element $x_{m+2}$ is not regular on $R/(J, x_{2m-1},\ldots,x_{m+3})$, and neither on $J_m[m-2]$ for that matter. And, in fact, we know that $J_m[m-2]$ has a whole batch of new linear syzygies forcing maximal linear rank.

(3) What are the odds against the first of the above conjectures? First, since one is in a homogeneous environment, one can just as well consider the sequence $\{x_{2m-1},\ldots, x_{m+3}\}$ in reverse order. At each step this way one is actually asking about the associated primes of $\bar{J}$, not those of $J_m[r]$ exactly.
Thus, it is urgent to compare both, a task that can be carried under severe hypothesis.
% (\cite{Huh}).
In any case, it would seem that the conjecture could use information about the associated primes of $J_m[r]$, a problem that has been  poorly accessed so far (Remark~\ref{Associated_primes}).
\end{Remark}

\begin{Question}\label{linear_info}\rm
	Assume that $m\geq 3$ and $0\leq r\leq m-2$.
  Is $J$ of linear type?
\end{Question}
One can check this up to small values of $m$ by computer calculation.
So far, the only proved case is $r=m-2$ (\cite{MAron}). Unfortunately, the argument there does not carry over to other values of $r$.
At this point it is not even clear that  $J$ satisfies property $(F_1)$, a front-runner feature of an ideal of linear type. In the generic case, $J$ is at least  a complete intersection locally at its unique minimal prime (\cite[Proof of Theorem 3.14]{MAron}). 

Proving affirmatively this question would give a short argument for both Theorem~\ref{non-homaloidal_main} (c) and Theorem~\ref{non_homaloidal} below.

\section{Special results in the generic case}\label{generic_case_revisited}

We now focus on the square generic Hankel matrix  $\mathcal{H}_m$ of order $m$ and on its determinant $f$. Let $R=k[\xx]=k[x_1,\ldots,x_{2m-1}]$ stand for the corresponding ground polynomial ring on the entries of $\mathcal{H}_m$.
Let $\ff$ stand for the set of derivatives of $f$.

\subsection{The background combinatorics}\label{combinatorics}

Recall that $\mathcal{H}_{r,s}$ denotes the $r\times s$ generic Hankel matrix.
When $r=s=:m$, one has $I_{m-1}(\mathcal{H}_{m,m}) =I_{m-1}(\mathcal{H}_{(m-1)\times(m+1)})$.
The advantage of this transfer is that not only one deals with the more pliable maximal minors (of the simplest case of a non-square generic Hankel matrix), but also the minimal number of generators becomes the predicted one ${{m+1}\choose 2}$.

For the current purpose, the maximal minor with columns $i_1<\cdots<i_{m-1}$ will be denoted $[i_1, \ldots, i_{m-1}]$, following a pretty much established notation.
By and large, one refers to some of the details developed in \cite[Section 3.3]{MAron}.
Pretty much as in the case of a generic matrix \cite[chapter 4]{B-V}, the set of maximal minors of the matrix $\mathcal{H}_{(m-1)\times(m+1)}$ is partially ordered by setting
$$[i_1, \ldots, i_{m-1}]\leq [j_1, \ldots, j_{m-1}]~~ \Leftrightarrow ~~ i_1\leq j_1, \ldots, i_{m-1}\leq j_{m-1}.$$

As an illustration, in the case of $m = 5$, one has

\smallskip

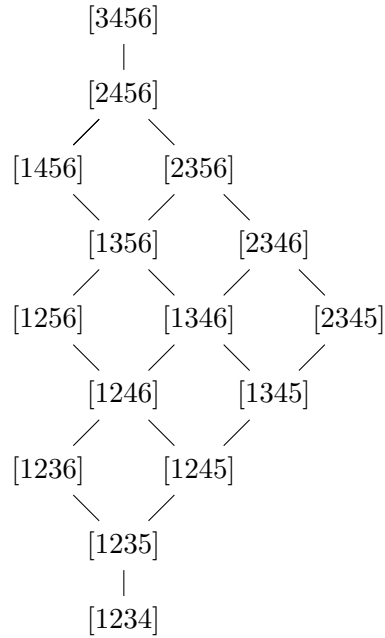
\hspace{2.2in}
\begin{tikzpicture}\label{diagram}
[scale=1.0,auto=center,every node/.style={circle,fill=blue!20}]
\node (n1) at (2,1) {[1234]};
\node (n2) at (2,2)  {[1235]};
\node (n3) at (1,3)  {[1236]};
\node (n4) at (3,3) {[1245]};
\node (n5) at (2,4)  {[1246]};
\node (n6) at (4,4)  {[1345]};
\node (n7) at (1,5)  {[1256]};
\node (n8) at (3,5)  {[1346]};
\node (n9) at (5,5)  {[2345]};
\node (n10) at (2,6)  {[1356]};
\node (n11) at (4,6)  {[2346]};
\node (n12) at (1,7)  {[1456]};
\node (n13) at (3,7)  {[2356]};
\node (n14) at (2,8)   {[2456]};
\node (n15) at (2,9)   {[3456]};
\foreach \from/\to in {n1/n2,n2/n3,n2/n4,n3/n5,n4/n5,n4/n6,n5/n7,n5/n8,n6/n8,n6/n9,n7/n10,n8/n10,n8/n11,n9/n11,n10/n12,n10/n13,n11/n13,n12/n14,n12/n14,n13/n14,n14/n15}
\draw (\from) -- (\to);
\end{tikzpicture}

where a link to the successive upper level  denotes $\leq$.

This poset has some remarkable properties reflected in the above diagram:

\begin{enumerate}
	\item[{\rm (1)}] Let ${n\choose 2}\leq l\leq {{m+2}\choose 2}-3$. The $l$th {\em level} of the poset is the subset of minors $[i_1,\ldots,i_{m-1}]$ such that $i_1+\cdots +i_{m-1}=l$ (i.e., the minors indexed by the ordered partitions of $l$)
	\item[{\rm (2)}]  A maximal minor in the $l$th level admits at most two comparable maximal minors on the $(l+1)$th level. To see this, note that a maximal minor $[i_1, \ldots, i_{m-1}]$ is alternatively designed by its complementary columns in the matrix. Thus, one can write 
	$$A=[i_1, \ldots, i_{m-1}]=[1,\ldots,i-1, \hat{i}, \ldots, j-1, \hat{j}, \ldots, m+1]$$
	where $\;\widehat{}\;$ denotes deletion. Therefore, in this notation there are at most two comparable minors with $A$ in the successive upper level, namely 
	$$[1,\ldots,\widehat{i-1}, i, \ldots, j-1, \hat{j}, \ldots, m+1] ~~\text{and} ~~[1,\ldots,i-1, \hat{i}, \ldots, \widehat{j-1}, j, \ldots, m+1].$$ 
\end{enumerate}

Thus far, the catalogue of the properties depends only on the size of the matrix, regardless as to what nature of specialization of the generic matrix one is considering.

For convenience, one sets  $y_{t,u}=x_j$, where $j=t+u-1, 1\leq j\leq 2m-1$. Thus, $\mathcal{H}_{m,m}=(y_{t,u})$.
Subtler properties of the maximal minors of $\mathcal{H}_{m-1,m+1}$ are as follows:
\begin{enumerate}
	\item[{\rm (3)}] (Level elements) The cofactors of $y_{t,u}$ and $y_{u,t}$ coincide and the maximal minors on the $l$th level of the poset are the distinct cofactors along the anti-diagonal $\{y_{t,u}| t+u=2m-l\}$.
	In particular, there is a ``central'' level with largest possible number of elements.
	\item[{\rm (4)}] (Partial derivatives) Let $f_j, 1\leq j\leq 2m-1,$ denote the partial derivative of $\det (\mathcal{H}_{m,m})$ with respect to $x_j$. Then $f_{2m-l}$ is a $k$-linear combination of the maximal minors in the $l$th level of the poset. 

	\item[{\rm (5)}] (Defining relations of the maximal minors)  In the generic case it is a classical result that the so-called {\em Pl\"ucker relations} generate the defining ideal of the maximal minors.
	The following result is slightly surprising:
	
	\begin{Proposition}\label{Hankel_is_Grassmann}%{\rm (char$k=0$)}
	Let $k$ be a field of characteristic zero and let $\mathcal{G}=(y_{u,v})_{\kern-5pt\tiny\begin{array}{c} 1\leq u\leq m-1\\ 1\leq v\leq m+1 \end{array}}$ and $\mathcal{H}=(x_{i+j-1})_{1\leq i\leq m-1, 1\leq j\leq m+1}$ denote, respectively,  the $(m-1)\times (m+1)$ generic matrix over $k$ and  the generic Hankel matrix of the same size over $k$.
	Then the respective special fiber algebras of the ideals $I_{m-1}(\mathcal{G})$ and $I_{m-1}(\mathcal{H})$ of maximal minors are isomorphic as graded $k$-algebras.
	\end{Proposition}
\demo
Consider the natural  $k$-algebra specialization map 
\begin{equation}\label{specialization_map}
k[y_{u,v}]\surjects k[x_1,\ldots,x_{2m-1}],\; y_{u,v}\mapsto x_{u+v-1}, \, 2\leq u+v\leq 2m
\end{equation}
of the corresponding ground polynomial rings.
Since each of the ideals in question is equigenerated, its fiber algebra is graded $k$-isomorphic to the respective $k$-subalgebra of the polynomial ring generated by the maximal minors, which one denotes here by $k[\boldsymbol{\Delta}(\mathcal{G})]$ and $k[\boldsymbol{\Delta}(\mathcal{H})]$, respectively.

Clearly, the map (\ref{specialization_map}) restricts to a $k$-algebra surjection \begin{equation}\label{fiber_algebra_surjection}
k[\boldsymbol{\Delta}(\mathcal{G})]\surjects k[\boldsymbol{\Delta}(\mathcal{H})].
\end{equation}

But both algebras have Krull dimension $2m-1$.

Indeed, the first because it is, up to grading normalization, graded isomorphic to the homogeneous coordinate ring of the Grassmannian $G(m-1,m+1)$, where the latter has the well-known dimension $2(m-1)=2m-2$.

The fact that also $\dim  k[\boldsymbol{\Delta}(\mathcal{H})]=2m-1$ as well follows from identifying the ideal of maximal minors $I_{m-1}(\mathcal{H})$ with the ideal of submaximal minors of the associated $m\times m$ Hankel matrix (see Proposition~\ref{GP}). After this is done the dimension comes out immediately off Proposition~\ref{max_anal_spread} or, alternatively, by the non-vanishing of the Hessian determinant of the square matrix as proved in \cite[Proposition 3.3.11]{Maral}, or yet by Corollary~\ref{Hessian_is_not_null}, whose proof is independent.

To conclude, the kernel of the surjection (\ref{fiber_algebra_surjection}) is a prime ideal of height zero in a domain, hence must be  the zero ideal.
\qed

\begin{Remark}\rm
For arbitrary size this turns out to be false, simply by having different dimensions, but perhaps more crucial as regards the present material, since minimal quadratic relations other than Pl\"ucker relations come in the picture as is explained in \cite{sagbi} and \cite{Conca1998}.
\end{Remark}
		\item[{\rm (6)}] (Shape of the Pl\"ucker relation)
	The Pl\"ucker relation in the case of a generic matrix of arbitrary size $p\times q$ has a somewhat cumbersome expression, with many indices floating around (see \cite[Lemma 4.4]{B-V}). Luckily, in the case of present interest, where $p=m-1,q=m+1$, assuming that the characteristic of the ground field is zero, the expression becomes simpler and these simpler expressions can be shown to generate the ideal of relations.
	This is based on the elementary observation to the effect that two maximal minors have in common at least $m-3$ indices. Thus, the intersection
	$$\{i_1,\ldots,i_{m-1}\}\cap (\{1,\ldots, m+1\}\setminus \{k_1,\ldots,k_{m-1}\})$$
	has at most $2$ elements.
	Consequently, the typical Pl\"ucker relation has at most $3$ terms, each a product of two minors, while the relevant ones correspond to the case where the above intersection has exactly $2$ elements, while the cases of $0$  and $1$ element yield, respectively,  an empty equation or to an identity.
	Thus, the shape can be described up to certain signs by assuming that, e.g., 
	$k_1=i_1,\ldots,k_{m-3}=i_{m-3}$ and the above intersection is the set $\{i_{m-2},i_{m-1}\}$, which affords the relation
	{\small
	 \begin{eqnarray*}\label{pluk}
	[i_1,\ldots,i_{m-1}]\kern-8pt&\cdot&\kern-8pt[k_1,\ldots,k_{m-1}]= [i_1,\ldots,i_{m-3},i_{m-2},i_{m-1}]\cdot [i_1,\ldots,i_{m-3},k_{m-2},k_{m-1}]\\
	&=&[i_1,\ldots,i_{m-3},\widehat{i_{m-2}},i_{m-1},k_{m-2}]\cdot [i_{m-2}, i_1,\ldots,i_{m-3},\widehat{k_{m-2}},k_{m-1}]\\
	&+& [i_1,\ldots,i_{m-3},i_{m-2},\widehat{i_{m-1}},k_{m-2}]\cdot [i_{m-1},  i_1,\ldots,i_{m-3},\widehat{k_{m-2}},k_{m-1}].
	\end{eqnarray*}
}	
\end{enumerate}

Up to signs and reordering of indices, this is the shape of a Pl\"ucker relation that intervenes in the proof of Theorem~\ref{is_the_radical}.

\begin{Remark}\rm
	As will be noted later, Pl\"ucker relations are not the right tool to deal with the degenerations $\mathcal{H}[r]$.
\end{Remark}

\subsection{The generic polar map is not homaloidal}\label{homaloidal}

According to common usage, the {\em polar map} of $f$ is the rational map $\pp^{2(m-2)}\dasharrow \pp^{2(m-2)}$ defined by the partial derivatives of $f$.
In a more classical terminology,  $f$ is a {\em homaloidal} polynomial if this map is birational, i.e., a Cremona map.

Dealing with these notions moves the emphasis  on to suitable $k$-subalgebras of $R$ rather than on its ideals.
Since $k[\ff]$ is a subalgebra of $R$ generated in degree $m-1$, the polar map is birational if and only if the homogeneous inclusion $k[\ff]\subset k[R_{m-1}]$ induces an isomorphism of the respective fields of fractions.

As a tool, we use yet another subalgebra, namely, let $\boldsymbol{\Delta}$ denote the set of $(m-1)$-minors  of  $\mathcal{H}_m$.
By Proposition~\ref{partials_vs_cofactors}, one has a homogeneous inclusion 
$k[\ff]\subset k[\boldsymbol{\Delta}]$. 

\begin{Theorem}\label{is_the_radical}
	$\sqrt{(\ff)k[\boldsymbol{\Delta}]}=(\boldsymbol{\Delta})k[\boldsymbol{\Delta}]$.
\end{Theorem}
\demo It clearly suffices to show that any minor $\Delta\in \boldsymbol{\Delta}$ has a power in $(\ff)k[\boldsymbol{\Delta}]$.
This fact has been essentially obtained in \cite[Proposition 3.13]{MAron} although the statement of the proposition there claimed less than was actually shown in the proof.

For convenience and precision we retrace the main parts of the argument in \cite[Proposition 3.13]{MAron} by pointing out the crucial re-editing.
Thus, one first lists $\ff=\{f_1,\ldots, f_{2m-1}\}$ and expresses each $f_j$ as a sum of maximal minors of  $\mathcal{H}_{m-1,m+1}$, as explained in property (4) above.

As no harm is done, one keeps denoting any of these maximal minors and their set by the same symbols $\Delta, \boldsymbol{\Delta}$, respectively.
The argument is then by descending induction on $j=1,\ldots,2m-1$ by showing that any maximal minor $\Delta$ which is a summand of $f_j$ has a power in $(\ff)k[\boldsymbol{\Delta}]$.

For $j=2m-1, 2m-2$ this is trivial since $f_{2m-1}, f_{2m-2}$ are themselves (signed) minors up to a coefficient.
As in the proof of \cite[Proposition 3.13]{MAron} and for later reasons, one displays the details of one additional explicit step in the induction -- although this step is really embodied in the general inductive step.
Namely, one deals now with the maximal minor $\Delta= [1,\ldots,n-3,n-1, n]$, a summand of the partial derivative $f_{2m-3}$.
For this, consider the following Pl\"ucker relation involving $\Delta$ and $\Delta'=[1,\ldots,n-2,n+1]$:
%{\small
\begin{equation}\label{plucker}
\Delta \Delta' =  1/2[1,\ldots,n-3,n-1,n+1]f_{2m-2}
-[1,\ldots,n-3,n,n+1] f_{2m-1},
\end{equation}
%}
where it will crucial that at least one element of each right-side factor belong to $J$.

On the other hand, by the above properties (3) and (4) together one knows that $f_{2m-3}$ is a $k$-linear (actually $\ZZ$-linear) combination of the minors $\Delta, \Delta'$.

Now,  take the obvious relation:
\begin{equation}\label{xequation}
\Delta^2 -1/3 \Delta ( \lambda\Delta + \Delta' ) + 1/\lambda \Delta \Delta'=0.
\end{equation}
where $\lambda$ is the suitable integer coefficient that appears in the latter $\ZZ$-linear combination.

From this follows immediately that $\Delta^2\in (\ff)k[\boldsymbol{\Delta}]$.

For the general inductive step the argument is similar, by showing that any $\Delta$ which appears as a summand of $f_j$, for $j\leq 2m-3$, satisfies an equation of integral dependence of the form
\begin{equation}\label{dependence_equation}
\Delta^2+a\Delta f_j+bg=0,
\end{equation}
for suitable $a,b\in k$, where $g$ is a (signed) sum of products of two maximal minors, one of which always appears as a summand of some $f_t$, with $t>j$. 
By the inductive hypothesis the latter minors have a power in $(\ff)k[\boldsymbol{\Delta}]$, and therefore so does every such $g$.
Finally, (\ref{dependence_equation}) implies that some power of $\Delta$ belongs to $(\ff)k[\boldsymbol{\Delta}]$ as well.
\qed

\begin{Remark}\rm
	In \cite[Proposition 3.13]{MAron} the stated result is that some power of $\Delta$ belongs to the {\em ideal}  $J=(\ff)$, but the argument clearly shows the stronger fact that  $\Delta$ belongs to $(\ff)k[\boldsymbol{\Delta}]$.
\end{Remark}

\begin{Corollary}\label{integrality}
	The ring inclusion  $k[\ff]\subset k[\boldsymbol{\Delta}]$ is an integral extension.
	In particular, $\dim k[\ff] =\dim k[\boldsymbol{\Delta}]$.
\end{Corollary}
\demo Fairly general, let $A=k[A_1]\subset B=k[B_1]$ denote a homogeneous inclusion of standard graded algebras over $k$. Then  $B_1\subset \sqrt{A_1B}$ implies that $B$ is finitely generated as an $A$-module.
In fact, letting $B_1^N\subset (A_1B)_N=A_1^N+A_1B^{N-1}$ for certain $N\geq 1$, then it is clear that $B$ is generated by $\{B_1,B_1^2,\ldots,B_1^{N-1}\}$ as an $A$-module.
\qed

\smallskip

One also retrieves the result of \cite[Proposition 3.3.11]{Maral}:

\begin{Corollary}\label{Hessian_is_not_null}
The Hessian determinant of $\mathcal{H}_m$ does not vanish. In particular, its partial derivatives are algebraically independent over $k$.
\end{Corollary}

The following theorem solves affirmatively some of the questions in \cite[Conjecture 3.16]{MAron}.

\begin{Theorem}\label{non-homaloidal_main}
	Let $\mathcal{H}_m$ denote the $m\times m \,(m\geq 3)$ generic Hankel matrix and let $k[\ff]\subset k[\boldsymbol{\Delta}]$ be as above.
	Then:
	\begin{enumerate}
		\item[{\rm (a)}] The fiber algebra $k[\boldsymbol{\Delta}]$ is a Gorenstein unique factorizations domain with regularity $m-2$.
		\item[{\rm (b)}] The extension $k[\ff]\subset k[\boldsymbol{\Delta}]$ is a Noether normalization.
			\item[{\rm (c)}] The polar map of $f=\det(\mathcal{H}_m)$ is not homaloidal.
		\item[{\rm (d)}] The ideal $J:=(\ff)\subset R$ is a reduction of the ideal of minors $I_{m-1}(\mathcal{H}_m)$ with reduction number $m-2$.
	\end{enumerate}
\end{Theorem}
\demo
(a) It is classical that the homogeneous coordinate ring of a Grassmannian is Gorenstein, as it follows from a result of \cite{Murthy}, by using that it is Cohen--Macaulay (\cite{Hochster}) and a unique factorization domain (\cite{Samuel}).

As for the regularity, one knows that the degree of the Hilbert series (as a rational function) of a Cohen--Macaulay standard graded $k$-algebra $C$ coincides with its $a$-invariant $a(C)$. On the other hand, under the same hypothesis, the Castelnuovo--Mumford regularity of $C$ coincides with the degree of the polynomial in the numerator of its Hilbert series. Since the denominator has degree $\dim C$, it follows that the regularity in the case of the Grassmannian (with the notation of  Proposition~\ref{Hankel_is_Grassmann}) is $2m-1-a(k[\boldsymbol{\Delta}(\mathcal{G})])$.
On the other hand, $a(k[\boldsymbol{\Delta}(\mathcal{G})])=-(m+1)$ (\cite[Corollary 1.4]{B-H_a}.

Therefore, the regularity of $k[\boldsymbol{\Delta}(\mathcal{G})]=2m-1-(m+1)=m-2$.
	
	To conclude, one uses Proposition~\ref{Hankel_is_Grassmann}.
	
(b) By Corollary~\ref{integrality} the extension $k[\ff]\subset k[\boldsymbol{\Delta}]$ is integral.
By \cite[Proposition 3.3.11]{Maral}, the Hessian of $f=\det(\mathcal{H}_m)$ does not vanish, hence $\ff$ is an algebraically independent set over $k$.
In other words, $k[\ff]$ is a polynomial ring over $k$.

(c) Since $k[\ff]$ is a polynomial ring it is integrally closed in its field of fractions. If the extension $k[\ff]\subset k[\boldsymbol{\Delta}]$ is birational then Corollary~\ref{integrality} implies that $k[\ff]= k[\boldsymbol{\Delta}]$.
But this is impossible since for $m\geq 3$ there is some cofactor which is not in the $k$-linear span of $k[\ff]$ (actually the $k$-linear span of the $(m-1)$-minors has dimension ${{m+1}\choose 2}>2m-1$ for $m\geq 3$).

(d) The first of the two statements is a well-known consequence of (b).
To see it, note that the $k$-algebra $k[\boldsymbol{\Delta}]$ is graded isomorphic to the special fiber of the ideal $I_{m-1}(\mathcal{H}_m)=(\boldsymbol{\Delta})\subset R$.
Then the result is explained in the proof of \cite[Theorem 1.77]{Wolmbook3}

For the reduction number, by \cite[Proposition 1.85]{Wolmbook3}, since $k[\boldsymbol{\Delta}]$ is Cohen--Macaulay, the reduction number of $J$ is the degree of the polynomial in the numerators of its Hilbert series.
But this is also the regularity by the facts in the proof of (a). Therefore, one gets the required value.
\qed

\begin{Remark}\rm
	A similar result to Theorem~\ref{non-homaloidal_main} (c) has been obtained by  Nivaldo Medeiros by a geometric argument.
\end{Remark}

\subsection{Elements of homaloidness for the degenerated polar map}

Recall, as mentioned before, that $f:=\det \mathcal H_m[r]$ is homaloidal for $r=m-2$ (\cite{CRS}).

The goal of this part is to give some elements towards the

\begin{Theorem}{\rm (Conjectured)}\label{non_homaloidal}	Assume that $0\leq r\leq m-2$. If  $f:=\det \mathcal H_m[r]$ is homaloidal then  $r=m-2$.
\end{Theorem}

The following collects a few elements detected in this regard.

\smallskip

{\sc Caveat 1.} The reason one cannot extend the method of proof of Theorem~\ref{non-homaloidal_main} (a) to the degenerated environment is that the extension $k[\ff]\subset k[\boldsymbol{\Delta}]$ is not integral anymore.
In fact, at least if $r\leq m-3$, then $P\not\subset \sqrt{J}$ since $R/J$ admits other minimal prime than $P$.
From the point of view of the actual argument of Theorem~\ref{non-homaloidal_main} (a), although one still has a poset of maximal minors, the existing Pl\"ucker relations do not necessarily come up with at least one factor belonging to $k[\ff]$.

\smallskip

 {\sc Caveat 2.} Since $k[\boldsymbol{\Delta}]$ has maximal dimension, for a degeneration $\mathcal{H}_m[r]$  one still has $\dim k[\ff]=\dim k[\boldsymbol{\Delta}]$ due to Theorem~\ref{Hessian_not_zero} -- of course, this equality is equivalent to the non-vanishing of the Hessian, so an a priori proof of the dimension equality would give a more elegant proof of the latter. 
Anyway, if one assumes homaloidness for $r\leq m-3$, then the corresponding field extension is trivial and since $k[\ff]$ is a polynomial ring, by \cite[Proposition 3.7 and Proposition 6.1 (a)]{ram2} one has the equality

$$e(k[\boldsymbol{\Delta}])=e(k[\ff], k[\boldsymbol{\Delta}])+1,$$
where $e(k[\ff], k[\boldsymbol{\Delta}])$ is the  relative multiplicity introduced in \cite{ram2}.
Note that the latter relative multiplicity can be computed as the ordinary multiplicity of the graded algebra
$$G/0:_G (\boldsymbol{\Delta})G,
$$
where $G$ stands for the associated graded ring of the ideal $(\ff)k[\boldsymbol{\Delta}]$ of the $k$-algebra $k[\boldsymbol{\Delta}]$.
In order to derive a contradiction one needs an a priori knowledge of a sufficiently large lower bound for $e(k[\boldsymbol{\Delta}])$, a presently unseen goal.

\smallskip

{\sc Caveat 3.} Picking up from a slightly different angle, still assuming that $f$ is homaloidal, let $\Delta\in \boldsymbol{\Delta}$ denote the cofactor of the $(m,m)$ entry,
the extension $A:=k[\ff]\subset B:=k[\ff, \Delta]$ is birational for even more reason.
Since $\Delta$ has same degree as the elements of $\ff$, it must be an element of the field of fractions $k(\ff)$ of the form $G(\ff)/F(\ff)$, with $F(\ff),G(\ff)\in k[\ff]$ homogeneous of degrees, say, $s-1,s$, respectively.
Since $k[\ff]$ is a polynomial ring over $k$, this forces a presentation of $B$ as a $k$-algebra over a polynomial ring $k[\mathbf{t}, u]$, with $\mathbf{t}=\{t_1,\ldots,t_{2m-r-1}\}$, as follows:

$$B\simeq k[\mathbf{t}, u]/(uF(\mathbf{t})+G(\mathbf{t})),$$
where $t_i\mapsto f_i, u\mapsto \Delta$, and $F,G$ are forms in $\mathbf{t}$ of degrees,  $s-1,s$, respectively.

The proof along these steps would then claim that this is impossible unless $r=m-2$.
In other words, one would show that such a presentation is only possible if $f_1=\partial f/\partial x_1$ is a pure power, in which case necessarily $f_1=x^{m-1}_{2m-r-1}=x^{m-1}_{m+1}$.

\smallskip
	
\smallskip 

{\sc Caveat 4:} It is equivalent to show that if $r< m-2$, the defining equation $e$ of $B$ associated to a presentation with set of generators $\{\ff,\Delta\}$ has $u$-degree $\geq 2$ (note that $e$ has no pure power term in $u$ since $u$ is not integral).

\smallskip

{\sc Caveat 5:} 
A computer calculation for the simplest situation ($m=4, r=1$) gives that $e$ has degree $15$ and $u$-degree $5$, well beyond what is needed. 

Note that either $F(\tt)$ or $G(\tt)$ must involve $t_1$ because the set $\{f_2,\ldots,f_{2m-r-1},\Delta\}$ is algebraically independent over $k$.
Suppose that $F(\tt)$ effectively involves $t_1$. Then by a $k$-linear change of variables fixing $t_1$ -- i.e., of the form $t_1\mapsto t_1, t_i\mapsto t_1+\alpha_it_i$, for $2\leq i\leq 2m-r-1$ -- one can assume that $F(\tt)$ has a non vanishing pure power term $\alpha t_1^{s-1}$. Then we'd have $\Delta f_1^{s-1}\in J^s$.
In this case it'd suffice to prove that for no exponent $l$ one has
$$\Delta\in J^l:f_1^{l-1}.$$

The feeling is that an early obstruction is the initial degree of this colon ideal which may turn out to be $>m-1$. Note that the colon ideal contains the ideal $J^l:J^{l-1}$, but one cannot derive anything from this as $J$ is seemingly Ratliff--Rush closed -- at least if one takes for granted that the associated graded ring $G$ of $J$ has positive depth given the expectation (Question~\ref{linear_info}) that $J$ is actually of linear type (not just analytically independent) and $G$ is Cohen--Macaulay.

As a slight confirmation, if $r=m-2$ is the case then a presentation equation as above can actually be written in the form $u F(\mathbf{t})+G(\mathbf{t})$ with $F=t_1^{e(B)-1}$, that is, $F(\ff)=f_1^{e(B)-1}=x_{m+1}^{(m-1)(e(B)-1)}$.

{\sc Caveat 6:} Early obstructions in the case $r<m-2$ are: (i) the ground polynomial ring $R$ has dimension at least $m+2$ (while the dimension is $m+1$ in the case $r=m-2$); (ii) $J\subset R$ -- and hence $(J,\Delta)\subset R$ too -- is an ideal of codimension $3$ (while $J$ has codimension $2$ and $(J,\Delta)$ has codimension $3$ in the case $r=m-2$).

\section{Appendix}

{\sc Proof of  the non-vanishing of the Hessian.}
The method consists in sufficiently ``degenerating'' the Hessian matrix to allow direct calculation with some of the submatrices to eventually get a non-vanishing expression. Namely, consider the ring endomorphism $\varphi$ of $R$ mapping any variable in $\vv:=\left\lbrace x_1, x_{m-r-1},x_{2m-r-1}\right\rbrace $ to itself and mapping any variable off $\vv$ to zero. We will show that by applying $\varphi$ to the entries of the Hessian matrix $H(f)$ the resulting matrix $H(f)(\vv)$ has non-vanishing determinant. %In other words, we will consider the entries of $H(f)$ modulo $\vv^c$ and will show that the determinant of the resulting matrix does not vanish.

For visualization we depict the matrix $\mathcal{H}_m[r]$ for arbitrary $r\leq m-3$:

\medskip

{\scriptsize
	$$\left( \begin{array}{ccccc|ccccc}
	x_1 & x_2 & \cdots & x_{m-r-2} & x_{m-r-1} & x_{m-r } & x_{m-r+1} & \cdots & x_{m-1} & x_m\\
	x_2 & x_3 & \cdots & x_{m-r-1} & x_{m-r} & x_{m-r+1} & x_{m-r+2} & \cdots & x_{m} & x_{m+1}\\
	\vdots  & \vdots  & \cdots &\vdots  &\vdots  & \vdots  & \vdots  & \cdots & \vdots  & \vdots \\
	x_{m-r-1} & x_{m-r } &\cdots  & x_{2m-2r-4}&x_{2m-2r-3} & x_{2m-2r-2} & x_{2m-2r-1}  & \cdots & x_{2m-r-3}  & x_{2m-r-2} \\
	\hline
	x_{m-r} & x_{m-r+1} &\cdots  & x_{2m-2r-3}&x_{2m-2r-2} & x_{2m-2r-1} & x_{2m-2r}  & \cdots & x_{2m-r-2}  & x_{2m-r-1}\\ 
	x_{m-r+1} & x_{m-r+2} &\cdots  & x_{2m-2r-2}&x_{2m-2r-1} & x_{2m-2r} & x_{2m-2r+1}  & \cdots & x_{2m-r-1}  & 0\\ 
	\vdots  & \vdots  & \cdots &\vdots  &\vdots  & \vdots  & \vdots  & \iddots & \vdots  & \vdots \\
	x_{m-1} & x_{m} &\cdots  & x_{2m-r-4}&x_{2m-r-3} & x_{2m-r-2} & x_{2m-r-1}  & \cdots & 0  & 0\\
	x_{m} & x_{m+1} &\cdots  & x_{2m-r-3}&x_{2m-r-2} & x_{2m-r-1} & 0  & \cdots & 0  & 0  
	\end{array}\right) $$
}

%In order to precisely locate an entry in the matrix it is convenient to reset the indices to double indices: 
%\begin{equation*}
%	z_{i,j} = \begin{cases}
%		x_{i+j-1},\quad \text{if $i+j\leq 2m-r$,}
%		\\
%		0, \quad \text{if $i+j>2m-r$}.
%	\end{cases}
%\end{equation*} 	

The goal is to isolate terms of the  partial derivatives of $f$ that have in their support a product of at least {\em two} variables off $\vv$, since such terms will produce at least one variable off $\vv$ in the entries of $H(f)$ and hence will vanish upon applying $\phi$.

The expression ``terms of degree at least $2$ off $\vv$'' will next appear recurrently in the sense just explained. In order to avoid tedious repeatition  we replace the expression by the letter $T$. 

Now, recall that the partial derivatives of $f$ are sums  of signed $(m-1)$-minors (see Proposition~~\ref{partials_vs_cofactors}). More precisely, 
for $k=1,\ldots,2m-r-1$, we have 
\begin{equation}\label{derivative_as_cofactor}
f_k=\sum_{i+j=k+1} M_{i,j},
\end{equation}
where $M_{i,j}$ is the (signed) cofactor of the $(i,j)$th entry. 

Let us pick up the shape of such a partial derivative of $f$ as we go through the various relevant intervals for the sum $i+j$.

Observe that  for  $i+j\leq m-r$, expanding the minor $M_{i,j}$ by the Laplace rule along its first $m-r-1$ rows yields 
$M_{i,j}= D_{i,j}x_{2m-r-1}^{r+1}+T,$ where $D_{i,j}$ is the cofactor of the $(i,j)$th entry of the submatrix

\begin{equation}\label{submat}
D=\left(
\begin{matrix}
x_1&x_2&\ldots &x_{m-r-1}\\
x_2&x_3&\ldots &x_{m-r}\\
\vdots &\vdots &\cdots &\vdots \\
x_{m-r-1}&x_{m-r}&\ldots &x_{2m-2r-3}\\
\end{matrix}
\right).
\end{equation}

\smallskip

\begin{Lemma}\label{lema1} Assume that $i+j\leq m-r$. Then:
	\begin{itemize}
		\item[\rm(a)] For $k+1=i+j<m-r$, one has 
		\begin{equation*}
		\varphi\left( \frac{\partial f_k}{\partial x_l}\right)=
		\begin{cases}
		\pm\, kx^{m-r-3}_{m-r-1}x_{2m-r-1}^{r+1},\, \text{if $l= 2m-2r-(k+1)-1$}\\
		0, \;\text{otherwise}
		\end{cases}
		\end{equation*}
		\item[\rm(b)] For $k+1=i+j=m-r$, one has
		%   $\varphi\left( \frac{\partial f_{m-r-1}}{\partial x_l}\right) $ :
		\begin{equation*}
		\varphi\left( \frac{\partial f_{m-r-1}}{\partial x_l}\right)=
		\begin{cases}
		\pm\, (m-r-1)(m-r-2)x_{m-r-1}^{m-r-3}x_{2m-r-1}^{r+1},\quad \text{if $l=m-r-1$} \medskip
		\\ %\vspace{5pt}
		\pm\, (m-r-1)(r+1)x_{m-r-1}^{m-r-2}x_{2m-r-1}^{r}, \quad \text{if $l= 2m-r-1$}\medskip\\
		\pm (m-r-3)\,x_1\,x_{m-r-1}^{m-r-4}\,x_{2m-r-1}^{r+1}, \quad \text{if $l=2m-2r-3$}\medskip\\
		0, \;\text{otherwise}
		\end{cases}.
		\end{equation*}

	\end{itemize}
\end{Lemma}
\demo 

(a) $i+j< m-r$

\smallskip

Expanding $D_{i,j}$ for $i+j<m-r$ yields: 
$$D_{i,j}=\pm\, z_{m-r-i,m-r-j}x_{m-r-1}^{m-r-3}+T= \pm\, x_{2m-2r-(i+j)-1}x_{m-r-1}^{m-r-3}+T.$$
Then (\ref{derivative_as_cofactor}) becomes
%\begin{equation}\label{derivative_as_cofactor2}
%f_k=\sum_{i+j=k+1} D_{i,j}x_{2m-r-1}^{r+1}+T,
%\end{equation}
%which then  gives, for $k+1=i+j< m-r$, 
\begin{eqnarray*}
	f_k&=& \sum_{i+j=k+1} D_{i,j}x_{2m-r-1}^{r+1}+T \\
	&=&\pm\, \sum_{i+j=k+1} x_{2m-2r-(i+j)-1}\,x_{m-r-1}^{m-r-3}\,x_{2m-r-1}^{r+1} +T\\
	&=& \pm\, k\cdot x_{2m-2r-(k+1)-1}\,x_{m-r-1}^{m-r-3}\,x_{2m-r-1}^{r+1}  +T.
\end{eqnarray*}

Clearly, $\varphi\left( \frac{\partial f_k}{\partial x_l}\right) $ vanish if $l\neq 2m-2r-(k+1)-1$. Moreover,  $$\varphi\left( \frac{\partial f_k}{\partial x_{2m-2r-(k+1)-1}} \right) =\pm kx^{m-r-3}_{m-r-1}x_{2m-r-1}^{r+1}$$

\medskip
%\begin{equation*}\phi\left(\frac{\partial^2f}{\partial x_l\partial x_k}\right)=
%\begin{cases}
%\pm\, k\cdot x_{m-r-1}^{m-r-3}x_{2m-r-1}^{r+1},\quad \text{if $l=2m-2r-(k+1)-1$}
%\medskip	\\ %\vspace{5pt}
%0, \quad \text{if $l\neq 2m-2r-(k+1)-1$}
%\end{cases}.
%\end{equation*} 

(b) For $i+j=m-r=k+1$ one has

$$D_{i,j}=
\left\lbrace \begin{array}{l}
\pm \, x_{m-r-1}^{m-r-2} +T, \;  \mbox{if $i=1$ or $j=1$} \medskip\\
\pm \, x_{m-r-1}^{m-r-2} \pm\, x_1 x_{2m-2r-3}\, x_{m-r-1}^{m-r-4}+T,\; \mbox{if $i\neq 1$ and $j\neq 1$ }\end{array}\right.. $$ 
Observe that, when $i\neq 1$ and $j\neq 1$ then we have necessarily $m-r\geq 4$. 

By (\ref{derivative_as_cofactor}) it obtains
\begin{eqnarray}\nonumber
f_{m-r-1}&=&\sum_{i+j=m-r} D_{i,j}x_{2m-r-1}^{r+1}+T=\pm (m-r-1)\,x_{m-r-1}^{m-r-2}\,x_{2m-r-1}^{r+1}\\ \nonumber
&\pm & (m-r-3)\,x_1\,x_{2m-2r-3}\,x_{m-r-1}^{m-r-4}\,x_{2m-r-1}^{r+1}  + T.
\end{eqnarray}

Taking derivative with respect to $x_l$ shows that  $\varphi\left( \frac{\partial f_{m-r-1}}{\partial x_{l}}\right) $ is  as  in the statement. 
\qed
\medskip
%\begin{equation*}
%\begin{cases}
%\pm\, (m-r-1)(m-r-2)x_{m-r-1}^{m-r-3}x_{2m-r-1}^{r+1},\quad \text{if $l=m-r-1$} \medskip
%\\ %\vspace{5pt}
%\pm\, (m-r-1)(r+1)x_{m-r-1}^{m-r-2}x_{2m-r-1}^{r}, \quad \text{if $l= 2m-r-1$}\medskip\\
%\pm (m-r-3)\,x_1\,x_{m-r-1}^{m-r-4}\,x_{2m-r-1}^{r+1}, \quad \text{if $l=2m-2r-3$}\medskip\\
%0,\quad \text{if $l\neq m-r-1,2m-r-1,2m-2r-3$} 
%\end{cases}.
%\end{equation*}  \qed

So far we have discussed the  first  $(m-r-1)$ columns  of $H(f)(\vv)$ and hence,  by symmetry, its first $(m-r-1)$ rows as well.
In particular,  the columns $m-r,\ldots, 2m-2r-3$ of $ H (f) (\vv) $ are partially obtained.
For the concluding argument on the entire shape of $H(f)(\vv)$ it will suffice to move all the way to the interval $2m-2r-1\leq i+j\leq  2m-r$, which will give the shape of columns $2m-2r-2, \ldots, 2m-r-1$ of  $ H (f) (\vv).$

\begin{Lemma}\label{lema2}  Assume that $2m-2r-1\leq i+j\leq  2m-r$. Then: \begin{itemize}
		\item[\rm(a)]  For $l=2m-2r-2,\ldots, 2m-r-2$ one has 	\begin{equation*}
		\phi\left(\frac{\partial f_l}{\partial x_{k}}\right)=
		\begin{cases}
		\pm\, 2{\bf q}=\pm 2 x_1x_{m-r-1}^{m-r-3}x_{2m-r-1}^r,\quad \text{if $k=4m-3r-l-4$,}
		\\ 0, \quad \text{if $k> 4m-3r-l-4$}.
		\end{cases}
		\end{equation*} 
		
		%If $k=4m-3r-l-4$, then  $\frac{\partial f_l}{\partial x_k}=$ module $\vv^C$ is $\pm2{\bf q}=\pm 2x_1x_{m-r-1}^{m-r-3}x_{2m-r-1}^r$ for $l=2m-2r-2,\ldots, 2m-r-2$. Moreover, $\frac{\partial f_l}{\partial x_k}=$ vanishes module $\vv^C$ if $k>4m-3r-l-4$.
		\item[\rm(b)] $\varphi\left( \frac{\partial f_{2m-r-1}}{\partial x_{2m-r-1}}\right)= \pm\,  (r+1)rx_{m-r-1}^{m-r-1}x_{2m-r-1}^{r-1}.
		$ 
		
	\end{itemize}
	
\end{Lemma}

\demo

(a) Let us look at the cofactors $M_{i,j}$ for $2m-2r-1\leq i+j< 2m-r$. Write $M_{i,j}=(-1)^{i+j}\det (C_{i,j})$, where $C_{i,j}$ is the submatrix of $\mathcal{H}_m(r)$ obtained by omitting its $i$th row and its $j$th column. 
If  $i> m-r-1$ and  $j> m-r-1$ then $C_{i,j}$ misses the entry $x_{2m-r-1}$ sitting on the  $(i,2m-r-i)$th and $(2m-r-j,j)$th slots $\mathcal{H}_m(r)$. 
Therefore,  in this case $C_{i,j}$ has only $(m-2)$ columns with some entry in $\vv$, and hence, $M_{i,j}$ cannot have any term supported on $\vv$; in addition, by the same token, any of its terms having degree $1$ in variables off $\vv$  involves necessarily the $m-r-1$ entries equal to $x_{m-r-1}$, the $r-1$ entries equal to $x_{2m-r-1}$ on the matrix $C_{i,j}$ and the entry of  $\mathcal{H}_m(r)$ in the $(2m-r-j,2m-r-i)$th slot.
But the latter is zero since  $(2m-r-j)+(2m-r-i)>2m-r$ when $2m-2r-1\leq i+j<  2m-r$.

Summing up, we have shown that, for $i> m-r-1$ and  $j> m-r-1$, the cofactor $M_{i,j}$ have neither terms supported on $\vv$ nor terms having degree $1$ in variables off $\vv$.

\smallskip

Thus, we are left with the next two possibilities: 

(i)  $i< m-r-1$ or $j<m-r-1$

By symmetry, it suffices to consider the case where $i< m-r-1$ and we do so.
Then  $C_{i,j}$ misses the entries $x_{m-r-1}$ and $x_{2m-r-1}$ in the $(i,m-r-i)$th and $(2m-r-j,j)$th slots of $\mathcal{H}_m(r)$, respectively. Clearly, $M_{i,j}$ does not have terms supported in $\vv$. Moreover, any term  of $M_{i,j}$ involving  $x_1$  cannot simultaneously involve variables of $\vv$ in slots $(m-r-1,1)$ and $(1,m-r-1)$ on $C_{i,j}$, and hence ought to have degree at least $2$ in the variables off $\vv$. Then one has  
\begin{eqnarray}\nonumber
%	M_{i,j}=\pm x_{m-r-1}^{m-r-2}x_{2m-r-1}^{r}a_{2m-r-j,m-r-i}+ T\\\nonumber
M_{i,j}=\pm x_{m-r-1}^{m-r-2}x_{2m-r-1}^{r}x_{3m-2r-(i+j)-1}+ T.
\end{eqnarray}

(ii) $i=m-r-1$  ou  $j=m-r-1$ 

Again, by symmetry it suffices to argue for the case where $i= m-r-1$.
A similar argument as above concerning slots $(m-r-1,1)$ and $(2m-r-j,j)$ of $\mathcal{H}_m(r)$ will do  and one gets 
$$\begin{matrix} M_{i,j}&=&\pm x_{m-r-1}^{m-r-2}x_{2m-r-1}^{r}x_{2m-r-j}\, \pm \;  x_1x_{m-r-1}^{m-r-3}x_{2m-r-1}^{r} x_{3m-2r-j-2} +  T\\
&=&\pm x_{m-r-1}^{m-r-2}x_{2m-r-1}^{r}x_{2m-r-j}\, \pm \;  x_1x_{m-r-1}^{m-r-3}x_{2m-r-1}^{r} x_{4m-3r-l-4} +  T
\end{matrix}$$
whereas $i+j=l$ e $i= m-r-1$.

A count of these cofactors gives for each $l=2m-2r-2,\ldots, 2m-r-2$:
\begin{equation*}
f_l=\pm 2 x_1x_{m-r-1}^{m-r-3}x_{2m-r-1}^rx_{4m-3r-l-4}\pm c_lx_{m-r-1}^{m-r-2}x_{2m-r-1}^rx_{3m-2r-l-2}+T.\nonumber
\end{equation*}	
for some $c_l\in k$. Therefore, for $l=2m-2r-2,\ldots, 2m-r-2$, one gets $$\varphi\left( \frac{\partial f_l}{\partial x_{4m-3r-l-4}}\right) = \pm\, 2{\bf q}=\pm 2 x_1x_{m-r-1}^{m-r-3}x_{2m-r-1}^r.$$
In addition, since  $3m-2r-l-2<4m-3r-l-4$ one sees that  $\varphi\left( \frac{\partial f_l}{\partial x_{k}}\right) =0$ if $k> 4m-3r-l-4$.
%\begin{equation*}
%\varphi\left( \frac{\partial f_l}{\partial x_{k}}\right) =
%\begin{cases}
%\pm\, 2{\bf q}=\pm 2 x_1x_{m-r-1}^{m-r-3}x_{2m-r-1}^r,\quad \text{if $k=4m-3r-l-4$,}
%\\ 0, \quad \text{if $k> 4m-3r-l-4$}.
%\end{cases}
%\end{equation*} 	

(b) One has  $k+1=i+j=2m-r$. Expanding along the first $m-r-1$ rows, one has $M_{i,j}= \det D\cdot x_{2m-r-1}^r+T$.
Expanding $\det D$ one gets  $$M_{i,j}=\pm x_1x_{m-r-1}^{m-r-3} x_{2m-r-1}^r x_{2m-2r-3}\pm x_{m-r-1}^{m-r-1}x_{2m-r-1}^r\ + T $$
and therefore, 
$$f_{2m-r-1}=\pm (r+1) x_1x_{m-r-1}^{m-r-3} x_{2m-r-1}^r x_{2m-2r-3}\pm (r+1)x_{m-r-1}^{m-r-1}x_{2m-r-1}^r\ + T. $$

Taking derivative with respect to $x_{2m-r-1}$ shows that
$$\hphantom{longlonglonglong}\varphi\left( \frac{\partial f_{2m-r-1}}{\partial x_{2m-r-1}}\right) = \pm\,  (r+1)rx_{m-r-1}^{m-r-1}x_{2m-r-1}^{r-1}.\hphantom{longlonglonglonglong}
\square $$

We now proceed to the proof of Theorem~\ref{Hessian_not_zero} proper.

Collecting the information gathered so far, we see that by applying $\phi$ to the entries of the Hessian matrix $H(f)$, one obtains a matrix in the form:

$$H(f)(\vv)=\left( \begin{array}{c|c}
A & B^t \\
\hline 
B& A' \end{array}\right). $$
Here $A$  and $B$ are matrices of sizes $(2m-2r-3)\times (2m-2r-3)$  and $(r+2)\times (2m-2r-3)$, respectively, and the stack $A\atop \overline{B}$ has the shape

{\scriptsize
	$$\left( \begin{array}{ccccccccc}
	0 & 0 &\ldots &	0 &  0 & 0 &  \ldots & 0 & \pm {\bf p}\\
	0 & 0 &\ldots &	 0&  0 & 0 &  \ldots   & \pm 2{\bf p}& 0\\
	\vdots  & \vdots  &\ldots &	  \vdots & \vdots & \vdots &  \iddots & \vdots& \vdots\\	 	
	0 & 0 &\ldots &	 0&  0 & \pm (m-r-2){\bf p} & \ldots & 0 &0\\	
	0 & 0 &\ldots &	 0&  \pm (m-r-1)(m-r-2){\bf p} & 0  &   \ldots  & 0& {\bf d} \\	
	0 & 0 &\ldots &	 \pm (m-r-2){\bf p}&  0 & *  &  \ldots & * & *\\	
	\vdots  & \vdots  &\iddots &	\vdots  & \vdots &  \vdots  & \ldots & \vdots & \vdots\\	
	0 & \pm 2{\bf p} &\ldots &	 0&  0 & *  &  \ldots & *&*\\	
	\pm {\bf p}   & 0&\ldots &	 0&  {\bf d}& *  & \ldots & *&*\\ [5pt]
	\hline\\ %[2pt]
	0 & 0 &\ldots &	 0&  0 & * &   \ldots  & *& * \\
	\vdots  & \vdots  &\iddots &	\vdots  & \vdots &  \vdots  & \ldots & \vdots & \vdots\\
	0 & 0 &\ldots &	 0&  0 & * &   \ldots  & *& * \\	
	0 & 0 &\ldots &	 0&  \pm (m-r-1)(r+1)x_{m-r-1}^{m-r-2}x_{2m-r-1}^{r}& * &   \ldots  & *& * 
	\end{array}\right) $$}
where ${\bf p}=x_{m-r-1}^{m-r-3} x_{2m-r-1}^{r+1}$ and ${\bf d}=\pm (m-r-3)\,x_1\,x_{m-r-1}^{m-r-4}\,x_{2m-r-1}^{r+1}$. This part follows from the Lemma (\ref{lema1}) 

As for the matrix $A'$, its shape follows from Lemma (\ref{lema2})
$$ A'=\left(\begin{array}{cccccc}
* & * & \ldots & \pm 2{\bf q}&0 & \\
\vdots  & \vdots & \iddots &\vdots  & \vdots  &\\
* & \pm 2{\bf q} & \ldots & 0& 0 & \\
\pm 2{\bf q} &  0& \ldots &0 & 0 &\\
0& 0 & \ldots &0 &  \pm (r+1)r x_{m-r-1}^{m-r-1}x_{2m-r-1}^{r-1}
\end{array} \right)  $$ 
with  ${\bf q}=x_1x_{m-r-1}^{m-r-3}x_{2m-2r-1}^r$.  

\medskip

Now expand the above determinant  along the first $2m-2r-3$  rows.  Note that the complementary minor to a $(2m-2r-3)$-minor of the first $2m-2r-3$ rows and avoiding the first $m-r-2$ columns vanishes as any of its columns is null.
At the other end,  the collection of non-vanishing minors of the first $(2m-2r-3)$ rows and involving the  first $m-r-2$ columns consists of $A$ itself and the following matrix that we will denote $X$:

{\tiny $$\left( \begin{array}{ccccccccc}
	0 & 0 &\ldots &	0 &  0 & 0 &0 &  \pm {\bf p}  & 0 \\
	0 & 0 &\ldots &	 0&   0  & 0 &\pm 2{\bf p} & 0 & 0\\
	\vdots  & \vdots  &\cdots &	  \vdots & \vdots &  \iddots & \vdots& \vdots & \vdots \\	 	
	0 & 0 &\ldots &	 0&  \pm (m-r-2){\bf p} & \ldots & 0 &0 & 0\\	
	0 & 0 &\ldots &	 0&    0  &\ldots & 0&  &\pm (m-r-1)(r+1)x_{m-r-1}^{m-r-2}x_{2m-r-1}^r  \\	
	0 & 0 &\ldots &	 \pm (m-r-2){\bf p}& *  &  \ldots & * & * & * \\	
	\vdots  & \vdots  &\iddots &	\vdots  &   \vdots  & \ldots & \vdots & \vdots & \vdots \\	
	0 & \pm 2{\bf p} &\ldots &	 0&   *  &  \ldots & *&* & *\\	
	\pm {\bf p}   & 0&\ldots &	 0&  *  & \ldots & *&* &  *\end{array} 
	\right) $$}

\noindent obtained upon replacing the $(m-r-1)$th column of $A$ with the last column of $B^t$ (i.e., the transpose of the last row of $B$).
Their complementary matrices are, respectively, $A'$ and $$X'=\left(\begin{array}{cccccc}
* & * & * & \ldots & \pm 2{\bf q} \\
\vdots &\vdots  & \vdots & \iddots &\vdots   \\
* & * & \pm 2{\bf q} & \ldots & 0  \\
*&\pm 2{\bf q} &  0& \ldots &0 \\
(m-r-1)(r+1)x_{m-r-1}^{m-r-2}x_{2m-r-1}^{r}&0& 0 & \ldots &0 
\end{array} \right)  $$  
Thereof, we obtain $ \det H(f)(\vv)=\pm \det A\det A'\pm \det X\det X'$.
Expanding the various determinants in this expression gives 
\begin{eqnarray*}
	\det H(f)(\vv)&=&  
	2^{r+1}(r+1)(m-r-1)!(m-r-2)!{\bf p}^{2m-2r-4}{\bf q}^{r+1}\\
	&\cdot &\kern-8pt\left(\pm r(m-r-2){\bf p} x_{m-r-1}^{m-r-1}x_{2m-r-1}^{r-1}\pm (m-r-1)(r+1)x_{m-r-1}^{2m-2r-4}x_{2m-r-1}^{2r}  \right).
\end{eqnarray*}
The first factor above is a  term in ${\bf p}$ and ${\bf q}$, hence does not vanish.
The second factor is a sum of distinct terms, so does not vanish either.
Therefore, the expression  is nonzero. 	
\qed

\begin{Remark}\rm 
	It would seem like there might exist an easy argument for the non-vanishing of the Hessian determinant of an intermediate $\det {\mathcal H}_m[r]$ since it is ``squeezed'' between the extreme situations where $r=0$ and $r=m-2$, where we know the Hessian does not vanish.
	Unfortunately, we may need the specifics of the present setup as for arbitrary threads of degenerations some intermediate Hessian determinants may vanish or not (see, e.g., \cite{CRS}, also \cite{apocriphal}). 
\end{Remark}

%%%%%%%%%%%%%%%%%%%%%%%%% bibliografia %%%%%%%%%%%%%%%%%%%%%%%%%%%%

\noindent {\bf Addresses:}

\medskip

\noindent {\sc Rainelly Cunha}\\
%Campus Natal\\
Instituto Federal de Educa\c c\~ao, Ci\^encia e Tecnologia do Rio Grande do Norte\\
59015-000  Natal, RN, Brazil\\
{\em e-mail}: rainelly.cunha@ifrn.edu.br \\

\noindent {\sc Maral Mostafazadehfard}\\
Instituto de Matem\'atica, CT--Bloco C\\
%Centro de Tecnologia - Bloco C, Cidade Universit\'aria\\
Universidade Federal do Rio de Janeiro\\
21941-909 Rio de Janeiro, RJ, Brazil\\
{\em e-mail}: maral@im.ufrj.br\\

\noindent {\sc Zaqueu Ramos}\\
Departamento de Matem\'atica, CCET\\ 
Universidade Federal de Sergipe\\
49100-000 S\~ao Cristov\~ao, Sergipe, Brazil\\
{\em e-mail}: zaqueu@mat.ufs.br\\

\noindent {\sc Aron Simis}\\
  Departamento de Matem\'atica, CCEN\\ 
Universidade Federal de Pernambuco\\ 
50740-560 Recife, PE, Brazil\\
{\em e-mail}:  aron@dmat.ufpe.br

\end{document}